\documentclass[11pt]{article}
\usepackage[margin=2cm]{geometry}
\usepackage{booktabs}
\usepackage{longtable}
\usepackage[T1]{fontenc}
\usepackage{lmodern}
\usepackage[numbers,longnamesfirst]{natbib}
\usepackage{algorithm,setspace}
\usepackage{algpseudocode}
\usepackage{subfigure} 
\usepackage{amsmath, amssymb, amsthm}
\usepackage{mathrsfs}
\usepackage{stmaryrd}
\usepackage{color, graphicx, float}
\usepackage[breakable]{tcolorbox}
\usepackage{tikz}
\usetikzlibrary{positioning, quotes}

\newlength{\R}

\usepackage{etoolbox}
\makeatletter
\patchcmd{\@maketitle}{\LARGE \@title}{\LARGE\bfseries\@title}{}{}

\renewcommand{\@seccntformat}[1]{\csname the#1\endcsname.\quad}
\makeatother

\usepackage{hyperref}
\hypersetup{
    colorlinks=true,
    linkcolor=red,      
    citecolor=blue,     
    urlcolor=blue       
}

\makeatletter
\def\th@plain{%
	\thm@notefont{}
	\itshape 
}
\def\th@definition{%
	\thm@notefont{}
	\normalfont 
}

\makeatother
\DeclareMathOperator{\inte}{int}
\newtheorem{theorem}{Theorem}[section]
\newtheorem{lemma}[theorem]{Lemma}

\theoremstyle{definition}
\newtheorem{definition}[theorem]{Definition}
\theoremstyle{definition}
\newtheorem{example}[theorem]{Example}
\theoremstyle{definition}
\newtheorem{remark}[theorem]{Remark}
\theoremstyle{definition}

\renewcommand{\labelenumi}{\rm (\roman{enumi})}

\usepackage[shortlabels]{enumitem}

\allowdisplaybreaks

\begin{document}

\title{Weak convergence of projection algorithm with momentum terms and new step size rule for quasimonotone variational inequalities}

\author{
Gourav Kumar\footnotemark[1],
~
Santanu Soe\footnotemark[1],
~
and~ V. Vetrivel\footnotemark[1]
}

\date{\today}

\makeatletter
\renewcommand{\@fnsymbol}[1]{1} 
\makeatother

\maketitle

\footnotetext[1]{Indian Institute of Technology Madras, Department of Mathematics, Chennai 600036, Tamil Nadu, India.\\
E-mails: \href{mailto:ma24r001@smai.iitm.ac.in}{ma24r001@smai.iitm.ac.in},
\href{mailto:ma22d002@smail.iitm.ac.in}{ma22d002@smail.iitm.ac.in},
\href{mailto:vetri@iitm.ac.in}{vetri@iitm.ac.in}.}

\maketitle

\begin{abstract}\noindent
This article analyses the simple projection method proposed by Izuchukwu et al. \cite[Algorithm 3.2]{izuchukwu2023simple} for solving variational inequality problems by incorporating momentum terms. A new step size strategy is also introduced, in which the step size sequence increases after a finite number of iterations. Under the assumptions that the underlying operator is quasimonotone and Lipschitz continuous, we establish weak convergence of the proposed method. The effectiveness and efficiency of the algorithm are demonstrated through numerical experiments and are compared with existing approaches from the literature. Finally, we apply the proposed algorithm to a signal recovery problem.
\end{abstract}


\noindent{\bf Mathematics Subject Classification (MSC 2020):}
47J20, 
49J40, 
65K15, 
65Y20.

\section{Introduction}
Let $\mathscr{C}$ be a nonempty, closed, and convex subset of a real Hilbert space $\mathscr{H}$, and $\mathscr{A}: \mathscr{H} \to \mathscr{H}$ be a continuous operator.
The objective of a variational inequality problem (VIP) is to find $p^*\in\mathscr{H}$ such that 
\begin{equation}\label{eqn_defvi}
    \langle \mathscr{A}p^*,p-p^*\rangle\geq 0~\forall p\in \mathscr{C}. 
\end{equation}
Problem (\ref{eqn_defvi}) was initially investigated by Fichera \cite{fichera1963sul} in the context of Signorini’s contact problem and was later independently examined by Stampacchia \cite{stampacchia}. Numerous significant problems in engineering, mathematical programming, economics, mechanics, etc. can be formulated as in the form of (\ref{eqn_defvi}) (see \cite{baiocchi1984variational,kinderlehrer2000introduction}, for more details).
\par
Many algorithms have been developed to compute solutions to the VIP (\ref{eqn_defvi}) under the assumption that the operator $\mathscr{A}$ is monotone (or pseudomonotone) and Lipschitz continuous. A widely used approach for this is the extragradient method \cite{korpelevich1976extragradient}: $\xi_1\in\mathscr{C}, \lambda_k\in \left(0,\frac{1}{L}\right)$, and $L>0$
\begin{align}\label{intro_1}
    &\left\{
    \begin{aligned}
        &\tau_k= \mathscr{P}_\mathscr{C}(\xi_k - \lambda_k \mathscr{A} \xi_k), \\
        &\xi_{k+1}= \mathscr{P}_\mathscr{C}(\xi_k - \lambda_k \mathscr{A} \tau_k), \quad k \geq 1.
    \end{aligned}
    \right.
\end{align}
Each iteration of (\ref{intro_1}) involves two evaluations of the operator $\mathscr{A}$ and two projection evaluations onto ${\mathscr{C}}$. As a result, the computational cost of the method can be significant, particularly in cases where $\mathscr{A}$ and $\mathscr{P}_{\mathscr{C}}$ are challenging to compute. To improve upon (\ref{intro_1}), Popov \cite{popov1980modification}, proposed the following extragradient method: $\xi_1, \eta_1\in\mathscr{C}$ and $\lambda_k\in \left(0,\frac{1}{3L}\right]$,
\begin{align}\label{intro_2}
    &\left\{
    \begin{aligned}
        &\xi_{k+1}= \mathscr{P}_\mathscr{C}(\xi_k - \lambda_k \mathscr{A} \tau_k), \\
        &\tau_{k+1}= \mathscr{P}_{\mathscr{C}}(\xi_{k+1} - \lambda_k \mathscr{A} \tau_k), \quad k \geq 1.
    \end{aligned}
    \right.
\end{align}
Unlike (\ref{intro_1}), in (\ref{intro_2}), the operator $\mathscr{A}$ is evaluated only once per iteration. However, $\mathscr{P}_\mathscr{C}$ is computed twice in each iteration. To mitigate this limitation, Censor et al. \cite{censor2012extensions} introduced a subgradient extragradient method that reduces the number of $\mathscr{P}_{\mathscr{C}}$ computations in Algorithm (\ref{intro_2}) to one per iteration: $\xi_1\in\mathscr{H}, \lambda_k\in\left(0,\frac{1}{L}\right)$ and $L>0$,
\begin{align}\label{intro_3}
    &\left\{
    \begin{aligned}
        &\tau_k= \mathscr{P}_\mathscr{C}(\xi_k - \lambda_k \mathscr{A} \xi_k), \\
        & D_k = \{w\in\mathscr{H}:\langle \xi_k-\lambda_k\mathscr{A}\xi_k-\tau_k, w-\tau_k\rangle\leq 0\},\\
        &\xi_{k+1}= \mathscr{P}_{D_k}(\xi_{k} - \lambda_k \mathscr{A} \tau_k), \quad k \geq 1.
    \end{aligned}
    \right.
\end{align}
Readers interested in other methods exhibiting similar properties to (\ref{intro_3}) are referred to \cite{gibali2020new,mainge2016convergence,malitsky2015projected}. The convergence analysis of method (\ref{intro_1}) under the assumption of pseudomonotonicity has been investigated in \cite{noor2003extragradient,vuong2018weak}. Weak convergence results for (\ref{intro_2}) when $\mathscr{A}$ is pseudomonotone are presented in \cite{censor2012extensions,malitsky2014extragradient}. Similar results for (\ref{intro_3}) have been examined in \cite{censor2011subgradient,kraikaew2014strong}.
 \par
In all the aforementioned algorithms, the step size sequence is either constant or nonincreasing. Malitsky \cite{malitsky2020golden} used the adaptive step size rule and proposed the adaptive golden ratio algorithm for (\ref{eqn_defvi}), which iterates as follows:
$\tau_0, \tau_1, \lambda_0 >0, \phi\in(1,\Phi), \bar{\lambda}>0$. Set $\xi_0=\tau_1, \theta_0=1, \rho=\frac{1}{\phi}+\frac{1}{{\phi}^2}$, and $k=1$.
\begin{align}\label{intro_4}
    &\left\{
    \begin{aligned}
        &\xi_{k}= \frac{(\phi-1)\tau_{k}+\tau_{k-1}}{\phi} \\
        &\tau_{k+1}= \mathscr{P}_{\mathscr{C}}(\xi_{k} - \lambda_k \mathscr{A} \tau_k),
    \end{aligned}
    \right.
\end{align}
where $\theta_k=\frac{\lambda_k\phi}{\lambda_{k-1}}$ and $\lambda_k=\min\left\{\rho\lambda_{k-1},\frac{\phi\theta_{k-1}}{4\lambda_{k-1}}\frac{\lVert\tau_k-\tau_{k-1}\rVert^2}{\lVert \mathscr{A}\tau_k-\mathscr{A}\tau_{k-1}\rVert^2},\bar{\lambda}\right\}.$ Method (\ref{intro_4}) retains all the advantages of (\ref{intro_3}) while additionally incorporating an adaptive step size strategy. The use of adaptive step size rules has been further investigated by several authors (see \cite{mtam,thong2021explicit,yang2021self}, and the references therein).
\par
It is important to note that only a limited number of studies have examined the aforementioned methods under the assumption that $\mathscr{A}$ is quasimonotone, a condition weaker than pseudomonotonicity. One of the primary challenges in this setting arises from the fact that the convergence analysis applicable to pseudomonotone operators does not directly extend to the quasimonotone case. For example, when $\mathscr{A}$ is quasimonotone, the dual variational inequality associated with (\ref{eqn_defvi}) (see (\ref{minty})) is no longer equivalent to (\ref{eqn_defvi}). Ye and He \citep[Algorithm 2.1]{ye2015double} proposed a double projection method and established its convergence to a solution of (\ref{eqn_defvi}) under the assumption that $\mathscr{A}$ is merely continuous in a finite-dimensional setting. Similar convergence results have also been demonstrated in \cite{tang2018strong,wang2020new}. Liu and Yang \cite{liu2020weak} established the weak convergence of the forward-backward-forward method to a solution of (\ref{eqn_defvi}) under the conditions that
$\mathscr{A}$ is quasimonotone, Lipschitz continuous, and sequentially weakly continuous in an infinite-dimensional Hilbert space. The aforementioned algorithms require projecting a vector onto the intersection of the feasible set 
$\mathscr{C}$ and $k+1$ half-spaces. However, as 
$n$ grows, this projection process becomes increasingly computationally demanding. To address this limitation, Izuchukwu et al. \citep[Algorithm 3.2]{izuchukwu2023simple} proposed a simple projection method, which iterates as follows.
\begin{algorithm}[H]
\caption{A Simple Projection Method for (\ref{eqn_defvi})}\label{algo_01} 
\begin{algorithmic}[1] 
\Statex {Input}: Let $\eta_0, \eta_1 > 0$, $\beta \in \left(\gamma, \frac{1 - 2\gamma}{2} \right)$ with $\gamma \in \left(0, \frac{1}{4} \right)$, and choose a nonnegative real sequence $\{u_k\}$ such that  
$
\sum_{n=1}^{\infty} u_k< \infty.
$ 
For arbitrary $\xi_0, \xi_1 \in \mathbb{C}$, let the sequence $\{\xi_k\}$ be generated by
\Statex Compute the next iterations:
\begin{align*}
\xi_{k+1}&= P_{\mathscr{C}} \left( \xi_k - \left( (\eta_k + \eta_{k-1}) A \xi_k +\eta_{k-1} A \xi_{k-1} \right) \right), \quad k \geq 1
,~\text{where}\\\eta_{n+1}&=
\begin{cases} 
\min \left( \alpha \frac{\| \xi_k - \xi_{k+1} \|}{\| A \xi_k - A \xi_{k+1} \|}, \eta_k + u_k \right), & \text{if } A \xi_k \neq A \xi_{k+1}\\  
\eta_k + u_k, & \text{otherwise}.
\end{cases}
\end{align*}
\end{algorithmic}
\end{algorithm}\noindent
With similar advantages as of Algorithm \ref{algo_01}, Izuchukwu and Shehu \citep[Algorithm 3]{izuchukwu2024golden} proposed a golden ratio algorithm that incorporates a backward inertial term.\par
To advance the solution methodology for quasivariational inequality problems, we introduce a modified version of Algorithm \ref{algo_01} enhanced with two momentum terms. Additionally, a nondecreasing step size sequence is employed to improve convergence behaviour. Under standard assumptions, we establish the weak convergence of the proposed algorithm. To validate the theoretical results, numerical experiments are conducted, and the performance is benchmarked against existing methods. The results demonstrate that the proposed approach outperforms several established algorithms in terms of number of iterations and CPU time.\par
The structure of the paper is as follows. Section \ref{sec:2} presents essential definitions and preliminary lemmas necessary for the convergence analysis. In Section \ref{sec:3}, we present our proposed method and discuss its weak convergence. Section \ref{sec:4} provides numerical experiments to evaluate the performance of the proposed method and the concluding remarks are given in Section \ref{sec:5}.
%

\section{Preliminaries}\label{sec:2}
\noindent
In this paper, the notation $v_k \rightharpoonup p$ is used to indicate that the sequence $(v_k)$ converges weakly to the point $p$.
\begin{definition}
    The operator $\mathscr{A}:\mathscr{H}\to\mathscr{H}$ is said to be
    \begin{enumerate}[(i)]
        \item\label{pre_itm_1} $L$-Lipschitz continuous if there exists $L>0$ such that
        \[\lVert\mathscr{A}p-\mathscr{A}q\rVert\leq L\lVert p-q\rVert~\forall p, q\in\mathscr{H};\]
        \item\label{pre_itm_2} monotone if 
        \[\langle \mathscr{A}p-\mathscr{A}q, p-q\rangle\geq 0~\forall p,q\in\mathscr{H};\]
        \item\label{pre_itm_3} pseudomonotone if 
        \[\langle\mathscr{A}q,p-q\rangle\geq 0\implies\langle\mathscr{A}p,p-q\rangle\geq 0~\forall p,q\in\mathscr{H};\]
        \item\label{pre_itm_4} quasimonotone if
        \[\langle\mathscr{A}q,p-q\rangle>0\implies\langle\mathscr{A}p,p-q\rangle\geq 0~\forall p,q\in\mathscr{H}.\]
    \end{enumerate}
\end{definition}
\noindent
It is known that $(\ref{pre_itm_2})\implies(\ref{pre_itm_3})\implies(\ref{pre_itm_4})$. However, the converse is not true.
A closely related problem to (\ref{eqn_defvi}) is to find $p^*\in\mathscr{C}$ such that
\begin{equation}\label{minty}
    \langle\mathscr{A}p,p-p^*\rangle\geq 0~\forall p\in\mathscr{C}.
\end{equation}
The VIP (\ref{minty}) is commonly referred to as the Minty VIP. Let $\mathscr{S}$ be the solution set of (\ref{eqn_defvi}) and $\mathscr{S}_\mathscr{M}$ be the solution set of (\ref{minty}). Then, $\mathscr{S}_\mathscr{M}$ is a closed and convex subset of $\mathscr{C}$. By the convexity of $\mathscr{C}$ and the continuity of $\mathscr{A}$, it follows that $\mathscr{S}_\mathscr{M}\subset\mathscr{S}$. The following lemma establishes sufficient conditions ensuring that $\mathscr{S}_\mathscr{M}$ is nonempty. 
\begin{lemma}\emph{\citep[Proposition 2.1]{ye2015double}} If either
\begin{enumerate}[(i)]
    \item $\mathscr{A}$ is pseudomonotone on $\mathscr{C}$ and $\mathscr{S}$ is nonempty,
    \item $\mathscr{A}$ corresponds to the gradient of a differentiable quasiconvex function $\mathscr{B}$,  which is defined over an open set $\mathscr{K}$ containing $\mathscr{C}$ and attains its minimum on $\mathscr{C}$,
    \item $\mathscr{A}$ is quasimonotone, nonzero and bounded on $\mathscr{C}$,
    \item $\mathscr{A}\neq 0$ and quasimonotone on $\mathscr{C}$ and there is $t>0$ for which $\lVert x\rVert\geq t$ for any $x\in\mathscr{C}$, there is $y\in\mathscr{C}$ satisfying $\lVert y\rVert\leq t$ and $\langle\mathscr{A}x,y-x\rangle\leq 0$,
    \item $\mathscr{A}$ is quasimonotone on $\mathscr{C}$, $\inte\mathscr{C}\neq\emptyset$, and there is $p^*\in\mathscr{S}$ for which $\mathscr{A}p^*\neq 0$,
\end{enumerate}
then $\mathscr{S}_\mathscr{M}\neq\emptyset.$
\end{lemma}
\noindent
The operator $\mathscr{P}_{\mathscr{C}}$ is defined on $\mathscr{H}$ with range in $\mathscr{C}$ and assigned to each $p\in\mathscr{H}$ the unique point in $\mathscr{C}$ that satisfies
\[\lVert p-\mathscr{P}_{\mathscr{C}}p\rVert=\inf_{q\in\mathscr{C}}\lVert p-q\rVert.\]
This operator is known as the metric projection. Moreover, for all $q\in\mathscr{C}$ the following inequality holds:
\begin{equation}\label{pre_2.2}
    \langle p-\mathscr{P}_{\mathscr{C}}p, q-\mathscr{P}_{\mathscr{C}}p\rangle\leq 0.
\end{equation}
The following lemmas are utilized in the subsequent discussion.
\begin{lemma}\label{pre_lm_1}\emph{\citep[Lemma 2.2]{izuchukwu2023simple}}
    Let $p, q\in\mathscr{H}$. Then, 
    \[2\langle p, q\rangle=\lVert p\rVert^2 +\lVert q\rVert^2-\lVert p-q\rVert^2=\lVert p+q\rVert^2-\lVert p\rVert^2-\lVert q\rVert^2.\]
\end{lemma}
\begin{lemma}\emph{\citep[Lemma 2.3]{hoai2024new}}\label{lm_pre_1}
    Let $(\xi_k)$ and $(\eta_k)$ be two sequences of nonnegative numbers satisfying:
    \[\xi_{k+1}\leq \xi_k-\eta_k~\forall k\in\mathbb{N}.\]
    Then, $(\xi_k)$ is convergent and $\sum_{k=0}^{\infty} \eta_k<\infty.$
\end{lemma}
\begin{lemma}\emph{\citep[Lemma 2.39]{bauschke2011convex}}\label{pre_lm_4}
    Let $(\xi_k)$ be a sequence in $\mathscr{H}$, and let $\mathscr{M}$ be a nonempty subset of $\mathscr{H}$. Assume that the following conditions hold:
    \begin{enumerate}[(i)]
        \item for every $p \in \mathscr{M}$, $\lim_{k\to\infty} \lVert \xi_k - p \rVert$ exists;
        \item every weak cluster point of $(\xi_k)$ is an element of $\mathscr{M}$.
    \end{enumerate}
    Then, there exists some $\bar{p} \in \mathscr{M}$ such that the sequence $(\xi_k)$ converges weakly to $\bar{p}$.
\end{lemma}

\section{Main results}\label{sec:3}\noindent
In this section, we establish the weak convergence for the proposed algorithm. To proceed, we impose the following standard assumptions.
\allowdisplaybreaks
\begin{enumerate}[(a)]
 \item\label{assum_1} $\mathscr{S}_\mathscr{M}\neq\emptyset$,
    \item\label{assum_2} $\mathscr{A}$ is $L$-Lipschitz continuous on $\mathscr{C}$,
    \item\label{assum_3} $\text{whenever } (w_k) \subset \mathscr{C} \text{ and } w_k \rightharpoonup p^*, \text{ we have } \|p^*\| \leq \liminf_{k \to \infty} \|w_k\|$,
    \item\label{assum_4} $\mathscr{A}$ is quasimonotone on $\mathscr{C}$.
\end{enumerate}
A formal description of the proposed algorithm is given below.
\begin{algorithm}[htbp]
\caption{Weakly Convergent Projection Algorithm with Momentum}\label{algo_1} 
\begin{algorithmic}[1] 
\Statex \textbf{Initialization}: Given $\theta\geq 0, \lambda_0=\lambda_1>0$, $\sigma\in(0,1)$ such that $\sigma<\frac{1}{3(1+\theta)}$ and a sequence $\gamma_k$ such that $\gamma_k\geq 0~\forall k\geq 0$ and $\sum_{k=0}^{\infty}\gamma_k< \infty.$
Choose $v_0, v_1, u_1\in \mathscr{H}$, and set $k=1$.
\Statex\textbf{Step 1}: Compute
\begin{align}
{w_k}=&~\frac{1}{1+\theta}v_k +\frac{\theta}{1+\theta}u_k\\
v_{k+1}=&~\mathscr{P}_{\mathscr{C}}(w_k-\lambda_k\mathscr{A}v_k-\lambda_{k-1}(\mathscr{A}v_k-\mathscr{A}v_{k-1}))\label{main_alg_step}\\
u_{k+1}=&~\frac{1}{1+\theta}v_{k+1}+\frac{\theta}{1+\theta}u_k.
\end{align}
\Statex \textbf{Step 2}: If 
$
    \lVert \mathscr{A}v_k-\mathscr{A}v_{k+1}\rVert>\frac{\sigma}{\lambda_{k}}\lVert v_k-v_{k+1}\rVert,
$
then
\begin{equation}\label{eqn_roc_2}
    \lambda_{k+1}=\sigma\frac{\lVert v_k-v_{k+1}\rVert}{\lVert \mathscr{A}v_k-\mathscr{A}v_{k+1}\rVert}
\end{equation}
else 
\begin{equation}\label{eqn_roc_3}
    \lambda_{k+1}=(1+\gamma_{k})\lambda_{k}.
\end{equation}
Set $k:=k+1$ and return to \textbf{Step 1}.
\end{algorithmic}
\end{algorithm}\noindent
\begin{remark}
\begin{enumerate}[(i)]
    \item Algorithm \ref{algo_1} requires only one metric projection onto $\mathscr{C}$ and a single evaluation of $\mathscr{A}$ per iteration, maintaining the same computational framework as the Algorithm \ref{algo_01}. This structure may offer numerical advantages in terms of iteration count when compared to (\ref{intro_1}) and its variants, which involve two evaluations of $\mathscr{A}$ per iteration.
    \item Unlike the methods proposed in \cite{censor2012extensions,korpelevich1976extragradient,popov1980modification}, Algorithm \ref{algo_1} does not require an estimate of the Lipschitz constant of the operator $\mathscr{A}$ as an input parameter.
\end{enumerate}
\end{remark}
\noindent
Following a similar proof as in Lemma 3.1 of \cite{hoai2024new}, it can be seen that the step size sequence $\lambda_k$ used in Algorithm \ref{algo_1}, satisfies the following property.
\begin{lemma}\label{lm_mm3}
  Suppose $\mathscr{A}$ is $L$-Lipschitz continuous and $(\lambda_k)$ is the step size sequence produced by {Algorithm} \ref{algo_1}. Then,
  there exists $k_1\in\mathbb{N}$ such that
        \[\lambda_{k+1}\geq\lambda_k~\forall k\geq k_1.\]
\end{lemma}\noindent
\begin{remark}\label{rmrk_1}
   Let the sequence $(\gamma_k)$ and the parameter $\sigma$ be defined as in Algorithm \ref{algo_1}. Since the series $\sum_{{k=0}}^\infty \gamma_k$ is convergent, there exists $\bar{k} \in \mathbb{N}$ such that
\[
\gamma_{k-1} \leq \frac{\sigma}{2}, \quad \forall k \geq \bar{k}.
\]
\end{remark}\noindent
We now introduce two fundamental lemmas that are essential for establishing the weak convergence of the algorithm.
\begin{lemma}\label{lm_mm2} Let $(v_k)$ be generated by Algorithm \ref{algo_1} and $p^*\in\mathscr{S}_{\mathscr{M}}$, then there exists $\bar{k}\in\mathbb{N}$ such that
\begin{enumerate}[(i)]
    \item\label{itm_lm_1} $2\lambda_k\langle \mathscr{A}v_k-\mathscr{A}v_{k-1}, v_k-v_{k+1}\rangle\leq\frac{3\sigma}{2}\Big(\lVert v_k-v_{k-1}\rVert^2 + \lVert v_k-v_{k+1}\rVert^2\Big)~\forall k\geq \bar{k},~\text{and}$
    \item\label{itm_lm_2} $2\lambda_k\langle \mathscr{A}v_k-\mathscr{A}v_{k-1}, v_k-p^*\rangle\leq\frac{3\sigma}{2}\Big(\lVert v_k-v_{k-1}\rVert^2 + \lVert v_k-p^*\rVert^2\Big)~\forall k\geq \bar{k}.$
\end{enumerate}
\end{lemma}\noindent
\allowdisplaybreaks
{\it \textbf{Proof} 
\ref{itm_lm_1}}
    We consider the following two possible cases.
    \begin{enumerate}[{Case} 1.]
        \item\label{lm_cs_1}  $\lVert \mathscr{A}v_k-\mathscr{A}v_{k-1}\rVert>\frac{\sigma}{\lambda_{k-1}}\lVert v_k-v_{k-1}\rVert.$\\
    By applying the Cauchy-Schwarz inequality and utilizing relation (\ref{eqn_roc_2}), we obtain 
        \allowdisplaybreaks
       \begin{align*}
                2\lambda_k\langle \mathscr{A}v_k-\mathscr{A}v_{k-1},v_k-v_{k+1}\rangle& \leq 2\lambda_k\lVert \mathscr{A}v_k-\mathscr{A}v_{k-1}\rVert~ \lVert v_k-v_{k+1}\rVert\\&=2\sigma\lVert v_k-v_{k-1}\rVert ~\lVert v_k-v_{k+1}\rVert\\&\leq\sigma\Big(\lVert v_k-v_{k-1}\rVert^2 + \lVert v_k-v_{k+1}\rVert^2\Big).
        \end{align*}
        \item \label{lm_cs_2} $ \lVert \mathscr{A}(v_k)-\mathscr{A}v_{k-1}\rVert\leq\frac{\sigma}{\lambda_{k-1}}\lVert v_k-v_{k-1}\rVert$.\\
        By applying the Cauchy-Schwarz inequality, relation (\ref{eqn_roc_3}), followed by Remark~\ref{rmrk_1}, we obtain the following $\forall k\geq\bar{k}$,
       \begin{align*}
                2\lambda_k\langle \mathscr{A}v_k-\mathscr{A}v_{k-1},v_k-v_{k+1}\rangle&\leq 2\lambda_k\lVert \mathscr{A}v_k-\mathscr{A}v_{k-1}\rVert~ \lVert v_k-v_{k+1}\rVert\\ &\leq 2\lambda_k\frac{\sigma}{\lambda_{k-1}}\lVert v_k-v_{k-1}\rVert~\lVert v_k-v_{k+1}\rVert\\&=2\sigma(1+\gamma_{k-1})\lVert v_k-v_{k-1}\rVert~\lVert v_k-v_{k+1}\rVert\\& <\sigma(1+\gamma_{k-1})\Big(\lVert v_k-v_{k-1}\rVert^2+\lVert v_k-v_{k+1}\rVert^2\Big)\\&\leq\frac{3\sigma}{2}\Big(\lVert v_k-v_{k-1}\rVert^2+\lVert v_k-v_{k+1}\rVert^2\Big).
        \end{align*}
    \end{enumerate}
     Thus, from Case \ref{lm_cs_1} and Case \ref{lm_cs_2}, we have
\[2\lambda_k\langle \mathscr{A}v_k-\mathscr{A}v_{k-1}, v_k-v_{k+1}\rangle\leq\frac{3\sigma}{2}\Big(\lVert v_k-v_{k-1}\rVert^2 + \lVert v_k-v_{k+1}\rVert^2\Big)~\forall k\geq \bar{k}.\]
The proof of \ref{itm_lm_2} follows on similar lines.
\qed
\allowdisplaybreaks
\begin{lemma}\label{lm_1}
    Under assumptions \ref{assum_1}-\ref{assum_4}, the sequences $(v_k), (w_k)$ and $(u_k)$ generated by Algorithm \ref{algo_1} are bounded. 
\end{lemma}
\allowdisplaybreaks\noindent
{\it \textbf{Proof}}
  Let $p^*\in\mathscr{S_\mathscr{M}}$. Then, $p^*\in\mathscr{S}\subset\mathscr{C}$. Thus, by (\ref{main_alg_step}), (\ref{pre_2.2}), and Lemma \ref{pre_lm_1}, we get 
  \begin{align}\label{lm_prf_1}
      0&\leq 2\langle v_{k+1}-(w_k-\lambda_k\mathscr{A}v_k-\lambda_{k-1}(\mathscr{A}v_k-\mathscr{A}v_{k-1})), p^*-v_{k+1}\rangle\nonumber\\&=2\langle v_{k+1}-w_k,p^*-v_{k+1}\rangle +2\lambda_{k-1}\langle \mathscr{A}v_k-\mathscr{A}v_{k-1},p^*-v_{k+1}\rangle\nonumber\\&~~~~+2\lambda_k\langle \mathscr{A}v_k,p^*-v_{k+1}\rangle\nonumber\\&= \lVert w_k-p^*\rVert^2-\lVert v_{k+1}-w_k\rVert^2-\lVert v_{k+1}-p^*\rVert^2+2\lambda_k\langle \mathscr{A}v_k,p^*-v_{k+1}\rangle\nonumber\\&~~~~+2\lambda_{k-1}\langle \mathscr{A}v_k-\mathscr{A}v_{k-1},p^*-v_{k+1}\rangle. 
  \end{align}
  Since $v_{k+1}\in\mathscr{C}$ and $p^*\in\mathscr{S}_{\mathscr{M}}\subset\mathscr{S}$, we get from (\ref{minty}) that $\langle \mathscr{A}v_{k+1},v_{k+1}-p^*\rangle\geq 0~\forall k\geq 1.$ This implies that $\langle \mathscr{A}v_k,p^*-v_{k+1}\rangle\leq\langle \mathscr{A}v_k-\mathscr{A}v_{k+1},p^*-v_{k+1}\rangle~\forall k\geq 1.$ Thus, (\ref{lm_prf_1}) becomes
  \begin{align}\label{lem_prf_2}
      \lVert v_{k+1}-p^*\rVert^2&\leq \lVert w_k-p^*\rVert^2-\lVert v_{k+1}-w_k\rVert^2+2\lambda_k\langle \mathscr{A}v_k-\mathscr{A}v_{k+1},p^*-v_{k+1}\rangle\nonumber\\&~~~~+2\lambda_{k-1}\langle \mathscr{A}v_k-\mathscr{A}v_{k-1},p^*-v_k\rangle+2\lambda_{k-1}\langle\mathscr{A}v_k-\mathscr{A}v_{k-1},v_k-v_{k+1}\rangle.
  \end{align}
  By using Lemmas \ref{lm_mm3} and \ref{lm_mm2}\ref{itm_lm_1}, $\forall k\geq k^*:=\max\{k_1,\bar{k}\}$, (\ref{lem_prf_2}) gives
  \begin{align*}
      \lVert v_{k+1}-p^*\rVert^2&\leq\lVert w_k-p^*\rVert^2-\lVert v_{k+1}-w_k\rVert^2 +2\lambda_k\langle \mathscr{A}v_k-\mathscr{A}v_{k+1},p^*-v_{k+1}\rangle\\&~~~~+2\lambda_{k-1}\langle \mathscr{A}v_k-\mathscr{A}v_{k-1},p^*-v_k\rangle +\frac{3\sigma}{2}\Big(\lVert v_k-v_{k-1}\rVert^2+\lVert v_k-v_{k+1}\rVert^2\Big).
  \end{align*}
  Thus, $\forall k\geq k^*$, we get
  \begin{align}\label{lem_prf_3}
      &\lVert v_{k+1}-p^*\rVert^2+2\lambda_k\langle \mathscr{A}v_{k+1}-\mathscr{A}v_k,p^*-v_{k+1}\rangle-\frac{3\sigma}{2}\lVert v_{k+1}-v_k\rVert^2\nonumber\\&\leq\lVert w_k-p^*\rVert^2-\lVert w_k-v_{k+1}\rVert^2+2\lambda_{k-1}\langle \mathscr{A}v_k-\mathscr{A}v_{k-1},p^*-v_k\rangle +\frac{3\sigma}{2}\lVert v_k-v_{k-1}\rVert^2.
  \end{align}
  Note that
  \begin{align}\label{lem_prf_4}
      \lVert w_k-p^*\rVert^2&=\left\lVert \frac{1}{1+\theta}(v_k-p^*)+\frac{\theta}{1+\theta}(u_k-p^*)\right\rVert^2\nonumber\\&=\frac{1}{1+\theta}\lVert v_k-p^*\rVert^2+\frac{\theta}{1+\theta}\lVert u_k-p^*\rVert^2-\frac{\theta}{(1+\theta)^2}\lVert v_k-u_k\rVert^2,
  \end{align}
  and
  \begin{align}\label{lem_prf_5}
      \lVert u_{k+1}-p^*\rVert^2&=\left\lVert\frac{\theta}{1+\theta}(u_k-p^*)+\frac{1}{1+\theta}(v_{k+1}-p^*)\right\rVert^2\nonumber\\&=\frac{\theta}{1+\theta}\lVert u_k-p^*\rVert^2+\frac{1}{1+\theta}\lVert v_{k+1}-p^*\rVert^2-\frac{\theta}{(1+\theta)^2}\lVert v_{k+1}-u_k\rVert^2.
  \end{align}
  Using (\ref{lem_prf_4}) into (\ref{lem_prf_3}), for all $k\geq k^*$, we get
  \begin{align}\label{lem_prf_6}
      &\lVert v_{k+1}-p^*\rVert^2 +2\lambda_k\langle \mathscr{A}v_{k+1}-\mathscr{A}v_k, p^*-v_{k+1}\rangle-\frac{3\sigma}{2}\lVert v_{k+1}-v_k\rVert^2\nonumber\\&\leq \frac{1}{1+\theta}\lVert v_k-p^*\lVert^2 +\frac{\theta}{1+\theta}\lVert u_k-p^*\rVert^2-\frac{\theta}{(1+\theta)^2}\lVert v_k-u_k\rVert^2\nonumber\\&~~~~-\lVert w_k-v_{k+1}\rVert^2 +2\lambda_{k-1}\langle \mathscr{A}v_k-\mathscr{A}v_{k-1},p^*-v_k\rangle+\frac{3\sigma}{2}\lVert v_k-v_{k-1}\rVert^2.
  \end{align}
  By multiplying (\ref{lem_prf_5}) by $\theta$ and summing with (\ref{lem_prf_6}), we obtain the following for all $k\geq k^*$,
  \begin{align}
      &\lVert v_{k+1}-p^*\rVert^2+\theta\lVert u_{k+1}-p^*\rVert^2+2\lambda_k\langle \mathscr{A}v_{k+1}-\mathscr{A}v_k,p^*-v_{k+1}\rangle-\frac{3\sigma}{2}\lVert v_{k+1}-v_k\rVert^2\nonumber\\&\leq \frac{1}{1+\theta}\lVert v_k-p^*\rVert^2+\frac{\theta}{1+\theta}\lVert u_k-p^*\rVert^2-\frac{\theta}{(1+\theta)^2}\lVert v_k-u_k\rVert^2-\lVert v_{k+1}-w_k\rVert^2\nonumber\\&~~~~+\frac{\theta}{1+\theta}\lVert v_{k+1}-p^*\rVert^2+\frac{{\theta}^2}{1+\theta}\lVert u_k-p^*\rVert^2-\frac{{\theta}^2}{(1+\theta)^2}\lVert v_{k+1}-u_k\rVert^2\nonumber\\&~~~~+2\lambda_{k-1}\langle \mathscr{A}v_k-\mathscr{A}v_{k-1},p^*-v_k\rangle+\frac{3\sigma}{2}\lVert v_k-v_{k-1}\rVert^2.
  \end{align}
Therefore, $\forall k\geq k^*$, we have
\begin{align}\label{lem_prf_7}
    &\frac{1}{1+\theta}\lVert v_{k+1}-p^*\rVert^2+\theta\lVert u_{k+1}-p^*\rVert^2+\frac{{\theta}^2}{(1+\theta)^2}\lVert v_{k+1}-u_k\rVert^2\nonumber\\&~~~~+2\lambda_k\langle \mathscr{A}v_k-\mathscr{A}v_{k+1},v_{k+1}-p^*\rangle\nonumber\\&\leq \frac{1}{1+\theta}\lVert v_k-p^*\rVert^2+\theta\lVert u_k-p^*\rVert^2-\frac{\theta}{(1+\theta)^2}\lVert v_k-u_k\rVert^2-\lVert v_{k+1}-w_k\rVert^2\nonumber\\&~~~~+\frac{{\theta}^2}{(1+\theta)^2}\lVert v_k-u_{k-1}\rVert^2-\frac{{\theta}^2}{(1+\theta)^2}\lVert v_k-u_{k-1}\rVert^2\nonumber\\&~~~~+2\lambda_{k-1}\langle \mathscr{A}v_{k-1}-\mathscr{A}v_k,v_k-p^*\rangle+\frac{3\sigma}{2}\lVert v_k-v_{k-1}\rVert^2+\frac{3\sigma}{2}\lVert v_{k+1}-v_k\rVert^2.
\end{align} 
Also, 
\begin{align}\label{lem_prf_8}
    \lVert v_{k+1}-w_k\rVert^2&=\left\lVert v_{k+1}-\frac{1}{1+\theta}v_k-\frac{\theta}{1+\theta}u_k\right\rVert^2\nonumber\\&=\left\lVert v_{k+1}-v_k-\frac{\theta}{1+\theta}(u_k-v_k)\right\rVert^2\nonumber\\&= \lVert v_{k+1}-v_k\rVert^2-\frac{2\theta}{1+\theta}\langle v_{k+1}-v_k,u_k-v_k\rangle+\left(\frac{\theta}{1+\theta}\right)^2\lVert u_k-v_k\rVert^2\nonumber\\&\geq \lVert v_{k+1}-v_k\rVert^2-\frac{2\theta}{1+\theta}\lVert v_{k+1}-v_k\rVert~\lVert u_k-v_k\rVert +\left(\frac{\theta}{1+\theta}\right)^2\lVert u_k-v_k\rVert^2\nonumber\\&\geq \lVert v_{k+1}-v_k\rVert^2-\frac{\theta}{1+\theta}\Big(\lVert v_{k+1}-v_k\rVert^2+\lVert u_k-v_k\rVert^2\Big)\nonumber\\&~~~~+\left(\frac{\theta}{1+\theta}\right)^2\lVert u_k-v_k\rVert^2\nonumber\\&=\frac{1}{1+\theta}\lVert v_{k+1}-v_k\rVert^2+\left(\left(\frac{\theta}{1+\theta}\right)^2-\left(\frac{\theta}{1+\theta}\right)\right)\lVert u_k-v_k\rVert^2.
\end{align}
Using (\ref{lem_prf_8}) in (\ref{lem_prf_7}) gives
\begin{align}
    &\frac{1}{1+\theta}\lVert v_{k+1}-p^*\rVert^2+\theta\lVert u_{k+1}-p^*\rVert^2+\frac{{\theta}^2}{(1+\theta)^2}\lVert v_{k+1}-u_k\rVert^2\nonumber\\&~~~~+2\lambda_k\langle \mathscr{A}v_k-\mathscr{A}v_{k+1},v_{k+1}-p^*\rangle\nonumber\\&\leq \frac{1}{1+\theta}\lVert v_k-p^*\rVert^2+\theta\lVert u_k-p^*\rVert^2-\frac{\theta}{(1+\theta)^2}\lVert v_k-u_k\rVert^2-\frac{1}{1+\theta}\lVert v_{k+1}-v_k\rVert^2\nonumber\\&~~~~-\left(\left(\frac{\theta}{1+\theta}\right)^2-\left(\frac{\theta}{1+\theta}\right)\right)\lVert u_k-v_k\rVert^2+\frac{{\theta}^2}{(1+\theta)^2}\lVert v_k-u_{k-1}\rVert^2\nonumber\\&~~~~-\frac{{\theta}^2}{(1+\theta)^2}\lVert v_k-u_{k-1}\rVert^2+2\lambda_{k-1}\langle \mathscr{A}v_{k-1}-\mathscr{A}v_k,v_k-p^*\rangle+\frac{3\sigma}{2}\lVert v_k-v_{k-1}\rVert^2\nonumber\\&~~~~+\frac{3\sigma}{2}\lVert v_{k+1}-v_k\rVert^2.
\end{align}
Therefore, $\forall k\geq k^*$, we have
\begin{align}\label{lem_prf_9}
    &\frac{1}{1+\theta}\lVert v_{k+1}-p^*\rVert^2+\theta\lVert u_{k+1}-p^*\rVert^2+\frac{{\theta}^2}{(1+\theta)^2}\lVert v_{k+1}-u_k\rVert^2\nonumber\\&~~~~+2\lambda_k\langle \mathscr{A}v_k-\mathscr{A}v_{k+1},v_{k+1}-p^*\rangle+\frac{3\sigma}{2}\lVert v_{k+1}-v_k\rVert^2\nonumber\\&\leq \frac{1}{1+\theta}\lVert v_k-p^*\rVert^2+\theta\lVert u_k-p^*\rVert^2+\frac{{\theta}^2}{(1+\theta)^2}\lVert v_k-u_{k-1}\rVert^2\nonumber\\&~~~~+2\lambda_{k-1}\langle \mathscr{A}v_{k-1}-\mathscr{A}v_k,v_k-p^*\rangle+\frac{3\sigma}{2}\lVert v_k-v_{k-1}\rVert^2\nonumber\\&~~~~-\left(\frac{1}{1+\theta}-3\sigma\right)\lVert v_{k+1}-v_k\rVert^2-\frac{{\theta}^2}{(1+\theta)^2}\lVert v_k-u_{k-1}\rVert^2.
\end{align}
Now, let us define
\begin{align}
    a_k:&=\frac{1}{1+\theta}\lVert v_k-p^*\rVert^2+\theta\lVert u_k-p^*\rVert^2+\frac{{\theta}^2}{(1+\theta)^2}\lVert v_k-u_{k-1}\rVert^2\nonumber\\&~~~~+2\lambda_{k-1}\langle \mathscr{A}v_{k-1}-\mathscr{A}v_k,v_k-p^*\rangle+\frac{3\sigma}{2}\lVert v_k-v_{k-1}\rVert^2.
\end{align}
Thus, from (\ref{lem_prf_9}), we have
\begin{equation}
    a_{k+1}\leq a_k-\left(\frac{1}{1+\theta}-3\sigma\right)\lVert v_{k+1}-v_k\rVert^2-\frac{{\theta}^2}{(1+\theta)^2}\lVert v_k-u_{k-1}\rVert^2.
\end{equation}
Thus, $a_k$ is monotonically decreasing.
By using Lemmas \ref{lm_mm3} and \ref{lm_mm2}\ref{itm_lm_2}, observe that
\begin{align}\label{bdd}
    a_k&=\frac{1}{1+\theta}\lVert v_k-p^*\rVert^2+\theta\lVert u_k-p^*\rVert^2+\frac{{\theta}^2}{(1+\theta)^2}\lVert v_k-u_{k-1}\rVert^2\nonumber\\&~~~~-2\lambda_{k-1}\langle \mathscr{A}v_k-\mathscr{A}v_{k-1},v_k-p^*\rangle+\frac{3\sigma}{2}\lVert v_k-v_{k-1}\rVert^2\nonumber\\&\geq\frac{1}{1+\theta}\lVert v_k-p^*\rVert^2+\theta\lVert u_k-p^*\rVert^2+\frac{{\theta}^2}{(1+\theta)^2}\lVert v_k-u_{k-1}\rVert^2\nonumber\\&~~~~-\frac{3\sigma}{2}\Big(\lVert v_{k-1}-v_k\rVert^2+\lVert v_k-p^*\rVert^2\Big)+\frac{3\sigma}{2}\lVert v_{k-1}-v_k\rVert^2\nonumber\\&=\left(\frac{1}{1+\theta}-\frac{3\sigma}{2}\right)\lVert v_k-p^*\rVert^2+\theta\lVert u_k-p^*\rVert^2+\frac{{\theta}^2}{(1+\theta)^2}\lVert v_k-u_{k-1}\rVert^2.
\end{align}
Thus, $a_k\geq 0~\forall k\geq k^*$. Therefore, according to Lemma \ref{lm_pre_1}, $\lim_{k\to\infty}a_k$ exists, and hence $(a_k)$ is bounded, also
\begin{equation}\label{hlm_prf_1}
    \lim_{k\to\infty}\lVert v_{k+1}-v_k\rVert=0=\lim_{k\to\infty}\lVert v_k-u_{k-1}\rVert.
\end{equation}
This implies from (\ref{bdd}) that $(v_k)$ and $(u_k)$ are bounded. By the definition of $w_k$ in Algorithm \ref{algo_1}, we have that $(w_k)$ is also bounded.
\qed\\
We now present the following lemma, the proof of which is inspired by the approach of Izuchukwu et al. \citep[Lemma 4.2]{izuchukwu2023simple}.
\begin{lemma}\label{lm_help_1}
Let $(v_k)$ be a sequence generated by Algorithm \ref{algo_1} under the assumptions \ref{assum_1}-\ref{assum_4}. Suppose that $p^*$ is a weak cluster point of $(v_k)$. Then, at least one of the following holds:  
\begin{enumerate}[(i)]
    \item $p^* \in \mathscr{S}_{\mathscr{M}}$, or  
    \item $\mathscr{A}p^* = 0$.
\end{enumerate}
\end{lemma}
\allowdisplaybreaks\noindent
{\it \textbf{Proof}}~
By Lemma \ref{lm_1}, $(v_k)$ is bounded. Therefore, let $p^*$ be a weak cluster point of $(v_k)$. Then, we can choose a subsequence of $(v_k)$, denoted by $(v_{{k}_{i}})$ such that $(v_{{k}_{i}})\rightharpoonup p^*\in\mathscr{C}.$ The analysis proceeds by considering two distinct cases. 
\begin{enumerate}[{Case} 1.]
\item Suppose that $\lim\sup_{k\to\infty}\lVert \mathscr{A}v_{{k}_{i}}\rVert=0.$ Then, $\lim_{k\to\infty}\lVert\mathscr{A}v_{{k}_{i}}\rVert=\lim\inf_{k\to\infty}\lVert\mathscr{A}v_{{k}_{i}}\rVert=0.$ Hence, from \ref{assum_3}, we get
\begin{equation}
    0\leq \lVert\mathscr{A}p^*\rVert\leq \lim\inf_{k\to\infty}\lVert\mathscr{A}v_{{k}_{i}}\rVert=0,
\end{equation}
which implies that $\mathscr{A}p^*=0.$
\item Suppose that $\lim\sup_{k\to\infty}\lVert\mathscr{A}v_{{k}_{i}}\rVert> 0.$ Then, without loss of generality, we can choose a subsequence of $\mathscr{A}v_{{k}_{i}}$ still denoted by $\mathscr{A}v_{{k}_{i}}$ such that $\lim_{k\to\infty}\lVert\mathscr{A}v_{{k}_{i}}\rVert=M_1>0$. By (\ref{main_alg_step}) and (\ref{pre_2.2}),  we get the following for all $x\in\mathscr{C}$ 
  \begin{align}\label{lm_prff_1}
      0&\leq \langle v_{k+1}-(w_k-\lambda_k\mathscr{A}v_k-\lambda_{k-1}(\mathscr{A}v_k-\mathscr{A}v_{k-1})), x-v_{k+1}\rangle\nonumber\\&=\langle v_{k+1}-w_k,x-v_{k+1}\rangle +\lambda_{k-1}\langle \mathscr{A}v_k-\mathscr{A}v_{k-1},x-v_{k+1}\rangle\nonumber\\&~~~~+\lambda_k\langle \mathscr{A}v_k,x-v_{k+1}\rangle. 
  \end{align}
  Since $\mathscr{A}$ is Lipschitz continuous on $\mathscr{C}$, from (\ref{hlm_prf_1}), we get
  \begin{equation}\label{hprf_2}
      \lim_{k\to\infty}\lVert \mathscr{A}v_k-\mathscr{A}v_{k-1}\rVert=0.
  \end{equation}
  Also, from (\ref{hlm_prf_1}) and the definition of $w_k$, we have
  \begin{equation}\label{hprf_3}
      \lim_{k\to\infty}\lVert v_{k+1}-w_k\rVert=0.
  \end{equation}
  Using (\ref{hprf_2}) and (\ref{hprf_3}) in (\ref{lm_prff_1}), we get 
  \[0\leq \liminf_{k\to\infty}\langle \mathscr{A}v_{{k}_{i}}, x-v_{{k}_{i}} \rangle\leq \limsup_{k\to\infty}\langle \mathscr{A}v_{{k}_{i}}, x-v_{{k}_{i}} \rangle<\infty~\forall x\in\mathscr{C}.\]
\end{enumerate}
The remainder of the proof proceeds analogously to that of \citep[Lemma~4.2]{izuchukwu2023simple}, and is therefore omitted.
\qed
\begin{theorem}
    The sequences $(v_k)$, $(w_k)$ and $(u_k)$ generated by Algorithm \ref{algo_1} converge weakly to a point in $\mathscr{S}_{\mathscr{M}}\subset \mathscr{S}$, when the assumptions \ref{assum_1}-\ref{assum_4}, and $\mathscr{A}x\neq 0~\forall x\in\mathscr{C}$ are satisfied.
\end{theorem}\noindent
{\it \textbf{Proof}}
Let the set of weak cluster points of $(v_k)$ be designated as $\mathscr{W}(v_k)$. Our aim is to establish that $\mathscr{W}(v_k)\subset\mathscr{S}_{\mathscr{M}}$. Let $p^*\in\mathscr{W}(v_k)$. By Lemma \ref{lm_1}, we can find a subsequence $(v_{k_{j}})\subset (v_k)$ for which $v_{k_{j}} \rightharpoonup p^*$ as $j\to \infty$. By the weakly closedness of $\mathscr{C}$, we arrive at $p^*\in\mathscr{C}$. Again, since for all $x\in\mathscr{C}$, $\mathscr{A}x$ is nonzero, we see that $\mathscr{A}p^*$ is also nonzero. Utilizing the Lemma \ref{lm_help_1}, we conclude that $p^*\in\mathscr{S}_{\mathscr{M}}$. Thus, $\mathscr{W}(v_k)\subset \mathscr{S}_{\mathscr{M}}$. Take 
\[b_k:=-2\lambda_{k-1}\langle \mathscr{A}v_{k-1}-\mathscr{A}v_k, v_k-p^*\rangle-\frac{3\sigma}{2}\lVert v_k-v_{k-1}\rVert^2\] and
\[c_k:=\theta\lVert u_k-v_k\rVert^2+2\theta\langle u_k-v_k,v_k-p^*\rangle+\frac{{\theta}^2}{(1+\theta)^2}\lVert v_k-u_{k-1}\rVert^2.\]
Using (\ref{hlm_prf_1}), and the fact that $\mathscr{A}$ is Lipschitz continuous, we have
\begin{equation}\label{end_1}
    \mathscr{A}v_{k-1}-\mathscr{A}v_k\to 0,~\text{as}~ k\to\infty.
\end{equation}
Also,
\begin{equation}\label{end_2}
    (v_k)~\text{is bounded, and again from (\ref{hlm_prf_1}) it can be seen that}~u_k-v_k\to 0~\text{as}~k\to\infty.
\end{equation}
From (\ref{hlm_prf_1}), (\ref{end_1}) and (\ref{end_2}), we conclude that $b_k\to 0$ as $k\to\infty$ and $c_k\to 0$ as $k\to\infty.$
Observe that 
\begin{align*}
    c_k=&\theta\lVert u_k-v_k\rVert^2+2\theta\langle u_k-v_k,v_k-p^*\rangle+\frac{{\theta}^2}{(1+\theta)^2}\lVert v_k-u_{k-1}\rVert^2\\=&\theta\lVert u_k-p^*\rVert^2-\theta\lVert v_k-p^*\rVert^2+\frac{{\theta}^2}{(1+\theta)^2}\lVert v_k-u_{k-1}\rVert^2.
\end{align*}
Then, $a_k-c_k+b_k=\left(\frac{1}{1+\theta}+\theta\right)\lVert v_k-p^*\rVert^2$.
Since $\lim_{k\to\infty}a_k$ exists and $\left(\frac{1}{1+\theta}+\theta\right)>0$, we have that $\lim_{k\to\infty}\lVert v_k-p^*\rVert$ exists. Therefore, by Lemma \ref{pre_lm_4}, $(v_k)$ converges weakly to a point in $\mathscr{S}_{\mathscr{M}}$. 
Furthermore, in light of (\ref{hlm_prf_1}) and (\ref{hprf_3}), we have that $(u_k)$ and $(w_k)$ converges weakly to a point in $\mathscr{S}_{\mathscr{M}}$. Therefore, the proof is concluded.
\qed
\section{Numerical results}\label{sec:4}\noindent
To assess the effectiveness of Algorithm \ref{algo_1}, we compare it with established approaches from the literature, specifically the algorithm of Izuchukwu and Shehu \citep[Algorithm 3]{izuchukwu2024golden}, Izuchukwu et al. \citep[Algorithm 3.2]{izuchukwu2023simple}, and the algorithm of Liu and Yang \citep[Algorithm 3.3]{liu2020weak}. The computational experiments are conducted using Python 3 on a personal computer equipped with an Intel(R) Core(TM) i7-2600 processor running at 2.30GHz and 8 GB of RAM. The parameters\footnote{The parameters for each algorithm are selected in accordance with their respective references \cite{izuchukwu2023simple,liu2020weak,ye2015double}} used to execute each algorithm are as follows:
\begin{enumerate}
\renewcommand{\labelenumi}{}
\item Algorithm 3 of \cite{izuchukwu2024golden}: $\delta=0.9, \theta=-0.1, \lambda_0=2, \bar{\lambda}=4$, tolerance, $\text{TOL}=\max\{\lVert w_{k+1}-v_k\rVert^2, \lVert w_k-v_k\rVert^2\}.$
\item Algorithm 3.2 of \cite{izuchukwu2023simple}: $\alpha=0.26, \delta = 0.17, a_k = \frac{100}{(k+1)^{1.1}}, \gamma_0=0.1, \gamma_1=0.01, ~\text{TOL}=\lVert v_k-P_{C}(v_k-\gamma(2Av_k-Av_{k-1}))\rVert + \lVert v_k-v_{k-1}\rVert, \gamma\in\left(\delta, \frac{1-2\delta}{L}\right).$
    \item Algorithm 3.3 of \cite{liu2020weak}: $\alpha=0.5, p_k=\frac{100}{(k+1)^{1.1}}$, $\lambda_0=0.01$, $\mu= 0.7$, $\text{TOL}=\frac{\lVert w_k-v_k\rVert}{\min\{\gamma_k,1\}}$. 
    \item Algorithm \ref{algo_1} (Our Algorithm): $\lambda_0=\lambda_1=0.01, \sigma=0.4\left(\frac{1}{2+2\theta}\right), \theta=0.01, \gamma_k=\frac{100}{(k+1)^{1.1}}, u_1=v_1$. In this case the tolerance is the same as Algorithm 3.2.
\end{enumerate}
We apply the stopping criterion $\text{TOL}<\epsilon$, where $\epsilon$ denotes the prescribed tolerance. Furthermore, the initial values $v_0, v_1$ are chosen as follows.\\
\\
For Example \ref{exm_1}: Case 1: $v_0=0.1, v_1=0.9$; Case 2: $v_0=0.8, v_1=0.1$; Case 3: $v_0=0.1, v_1=0.5$, and Case 4: $v_0=-0.1, v_1=0.2$.\\
For Example \ref{exm_2}: Case 1: $v_0=(0.3,0.1), v_1=(0.1,0.5)$; Case 2: $v_0=(0.1,0.1), v_1=(0.1,0.7)$; Case 3: $v_0=(0.1,-0.5), v_1=(0.1,0.3)$, and Case 4: $v_0=(0.3,-0.7), v_1=(0.2,-0.5)$.\\
For example \ref{exm_3}: $m\in(50, 80, 100, 200)$ while $v_0$ and $v_1$ are randomly generated.\\
For example \ref{exm_4}: Case 1: $v_0=\frac{1}{3^k}, v_1=\frac{2^k}{3^k}$; Case 2: $v_0=\frac{1}{2^k}, v_1=\frac{1}{5^k}$; Case 3: $v_0=\frac{4^k}{5^k}, v_1=\frac{1}{2^k}$, and Case 4: $v_0=\frac{1}{8^k}, v_1=\frac{1}{7^k}$.\par
We now present the following test problems for numerical evaluation, which have also been considered in prior studies such as \cite{izuchukwu2024golden, izuchukwu2023simple, liu2020weak}.
\begin{example}\label{exm_1}
    Let $\mathscr{C}=[-1,1]$ and 
    \[\mathscr{A}{\mu}=\begin{cases}
        2\mu-1, &\text{if}~\mu>1,\\
        {\mu}^2, &\text{if}~\mu\in[-1,1],\\
        -2\mu-1, &\text{if}~\mu<-1.
    \end{cases}\]
    It can be seen that $\mathscr{A}$ is quasimonotone and Lipschitz continuous. Also, $\mathscr{S}_\mathscr{M}=\{-1\}$ and $\mathscr{S}=\{-1,0\}.$
\end{example}
\begin{figure}[h]
\centering
\mbox{\subfigure{\includegraphics[scale=0.33]{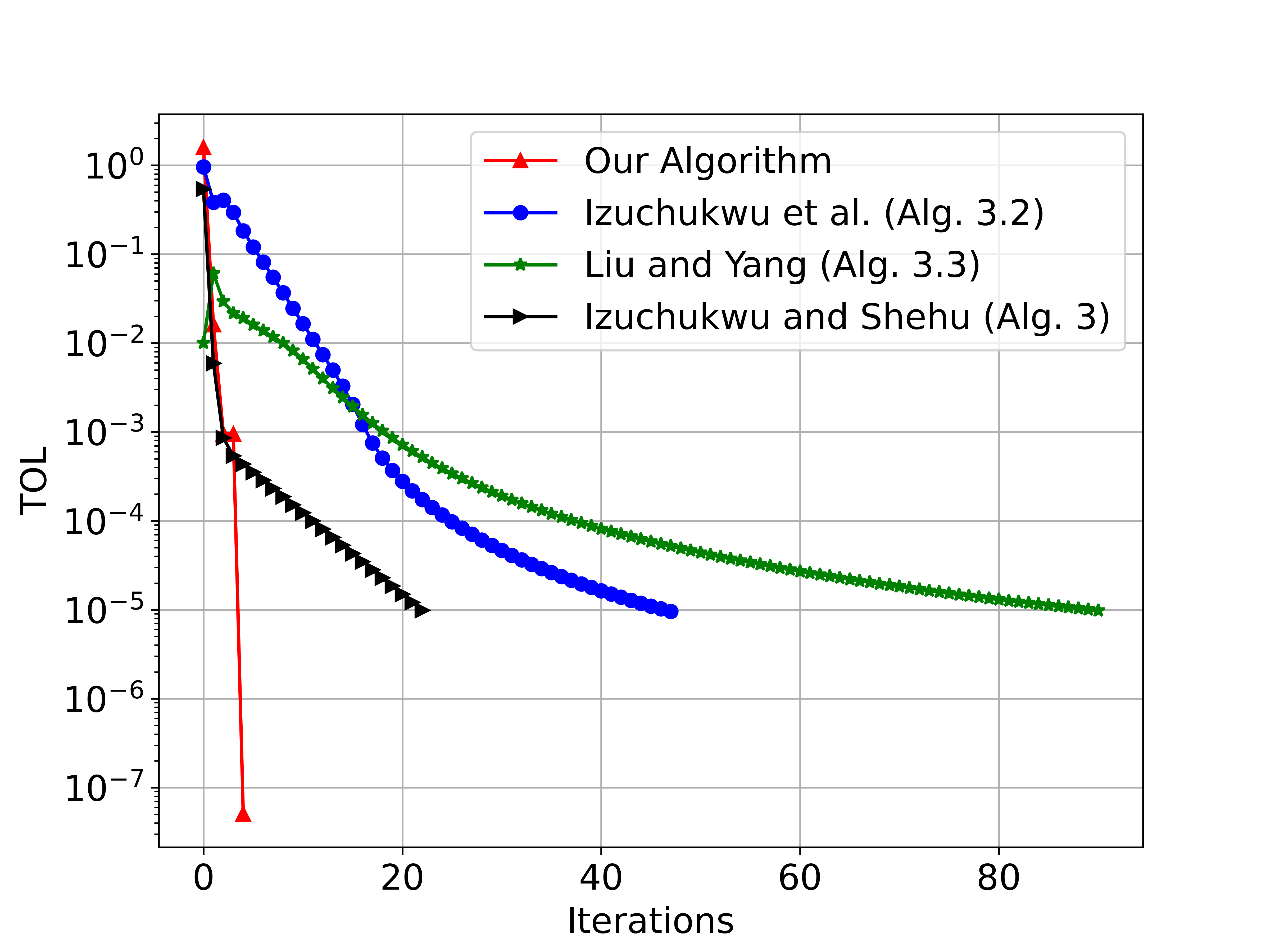} }\quad
\subfigure{{\includegraphics[scale=0.33]{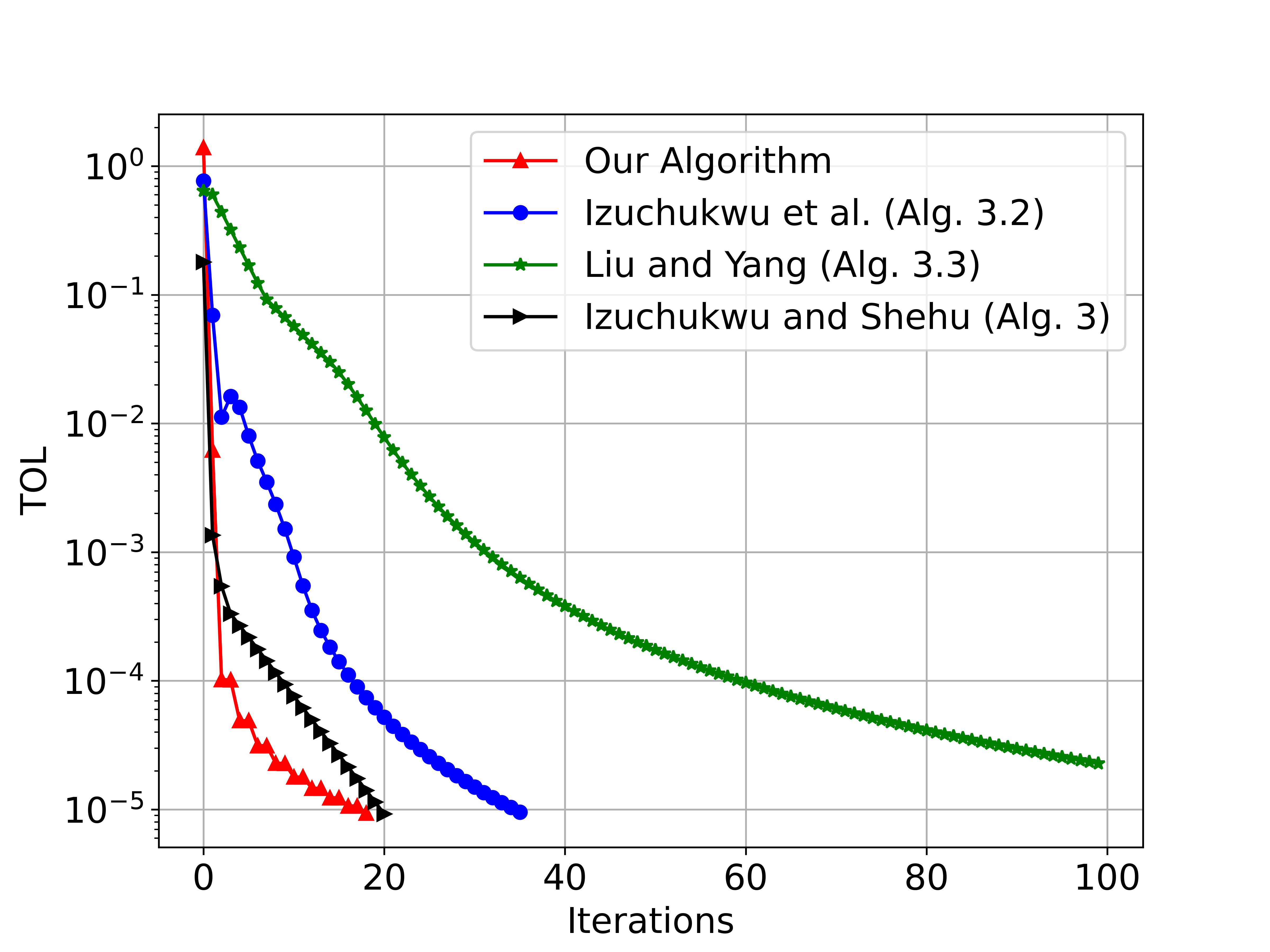}}}}\quad
\subfigure{{\includegraphics[scale=0.33]{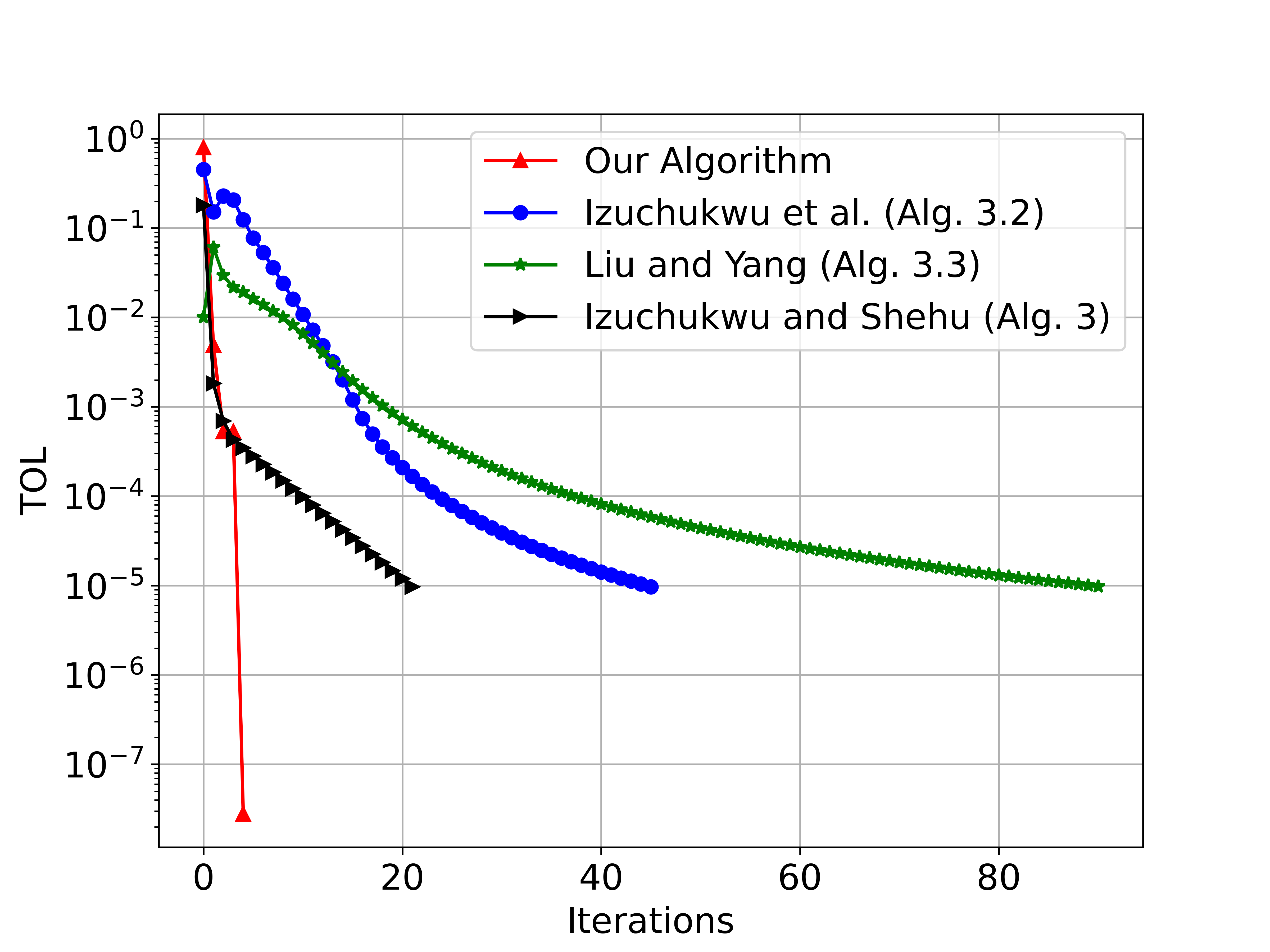}}}\quad
\subfigure{{\includegraphics[scale=0.33]{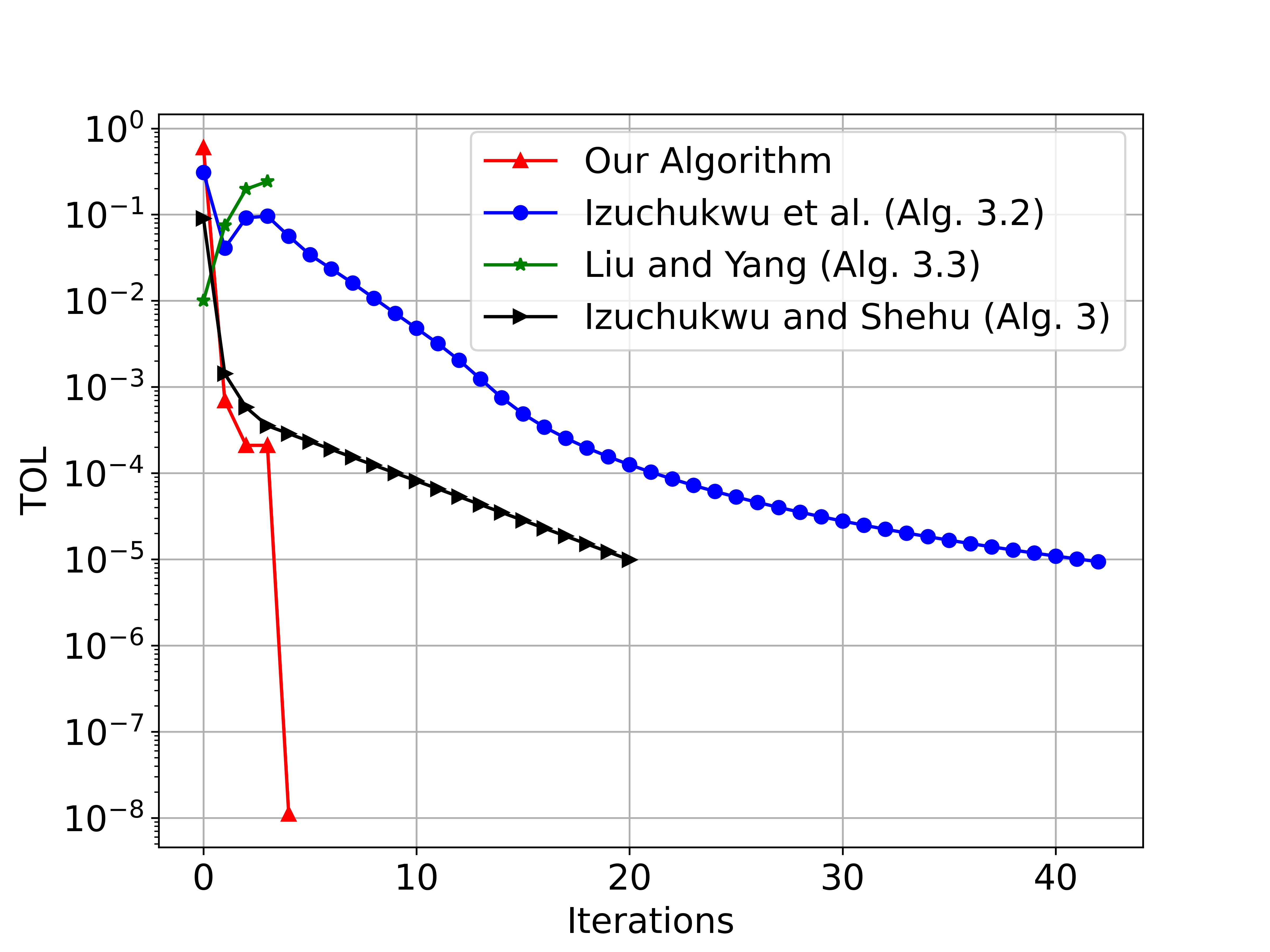}}}
\caption{Comparison of \text{TOL} for Example \ref{exm_1} with $\epsilon=10^{-5}$: Top Left: Case 1; Top Right: Case 2; Bottom left: Case 3; Bottom right: Case 4}
 \label{example_1}
\end{figure}
\begin{example}\label{exm_2}
    Let $\mathscr{C}=\{\mu\in{\mathbb{R}}^2:{\mu_1}^2+{\mu_2}^2\leq 1, 0\leq \mu_1\}$ and $\mathscr{A}(\mu_1,\mu_2)=(-\mu_1 e^{\mu_2},\mu_2)$. In this case, the operator $\mathscr{A}$ is not quasimonotone \citep[Problem 2]{liu2020weak}. Furthermore, it can be verified that $(1,0)\in\mathscr{S}_\mathscr{M}$ and $\mathscr{S}=\{(1,0), (0,0)\}.$
\end{example}
\begin{example}\label{exm_3}
    Let $\mathscr{C} = [0, 1]^m$ and $\mathscr{A}\mu = (h_1\mu, h_2\mu, \ldots, h_m\mu)$, where
\[
h_i\mu = \mu_{i-1}^2 + \mu_i^2 + \mu_{i-1}\mu_i + \mu_i\mu_{i+1} - 2\mu_{i-1} + 4\mu_i + \mu_{i+1} - 1,\quad i = 1, 2, \ldots, m,
\]
with boundary conditions $\mu_0 = \mu_{m+1} = 0$.
\end{example}
\begin{example}\label{exm_4}
   Next, we take an example in an infinite dimensional Hilbert space. Let $\mathscr{H}=\mathscr{L}_2:=\{x=(x_1,x_2, \ldots, x_i,\ldots):\sum_{i=1}^{\infty}|x_i|^2<\infty\}$. Take $\mathscr{C}=\{x\in\mathscr{L}_2:\lVert x\rVert\leq 3\}$ and $\mathscr{A}x:=(x_1e^{-{x_1}^2},0,0,\ldots), x=(x_1,x_2,x_3,\ldots)\in\mathscr{C}$. It has been demonstrated in \citep[Example 5.4]{izuchukwu2023simple} that the operator $\mathscr{A}$ satisfies the property of quasimonotonicity.
    \end{example}
\begin{table}[htbp]
\centering
\caption{Example \ref{exm_1}: Comparison of algorithms with $\epsilon = 10^{-5}$}
\begin{tabular}{lcccccccc}
\toprule
Algorithms & \multicolumn{2}{c}{Case 1} & \multicolumn{2}{c}{Case 2} & \multicolumn{2}{c}{Case 3} & \multicolumn{2}{c}{Case 4} \\
\cmidrule(r){2-3} \cmidrule(r){4-5} \cmidrule(r){6-7} \cmidrule(r){8-9}
 & CPU(s) & Iter & CPU(s) & Iter & CPU(s) & Iter & CPU(s) & Iter \\
\midrule
Our Algorithm & 0.0006 & 5 & 0.0011 & 19 & 0.0006 & 5 & 0.0004 & 5 \\
Algorithm 3 in \cite{izuchukwu2024golden}& 0.0027 & 23 & 0.0012 & 21 & 0.0012 & 22 & 0.0012 & 21 \\
Algorithm 3.2 in \cite{izuchukwu2023simple}& 0.0057 & 48 & 0.0028 & 36 & 0.0048 & 45 & 0.0034 & 43 \\
Algorithm 3.3 in \cite{liu2020weak}& 0.0160 & 91 & 0.0074 & 100 & 0.0115 & 91 & 0.0012 & 5\\
\bottomrule
\end{tabular}
\end{table}
%
%
%
%
%
%
\begin{figure}[h]
\centering
\subfigure{\includegraphics[scale=0.35]{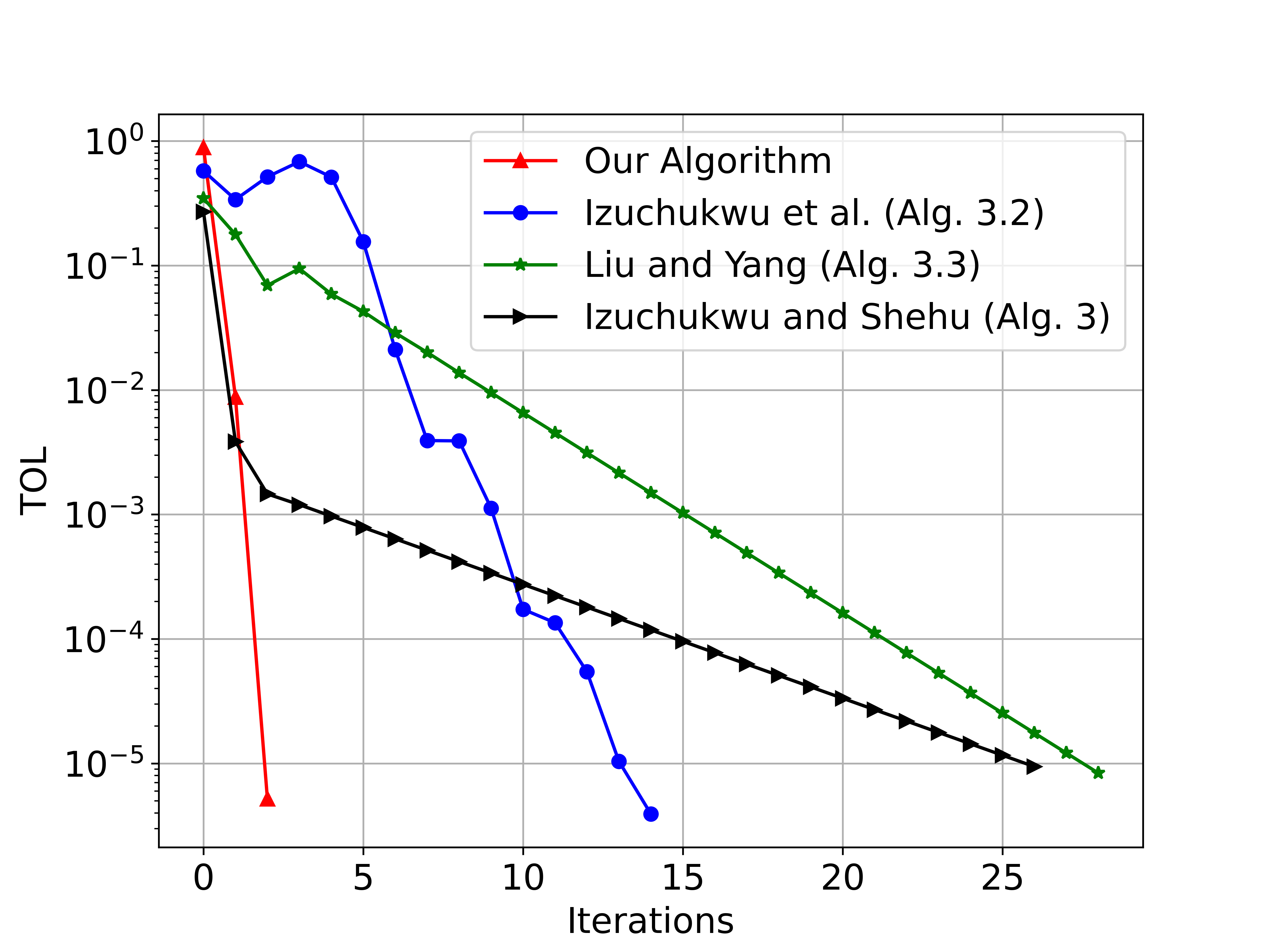} }\quad
\subfigure{{\includegraphics[scale=0.35]{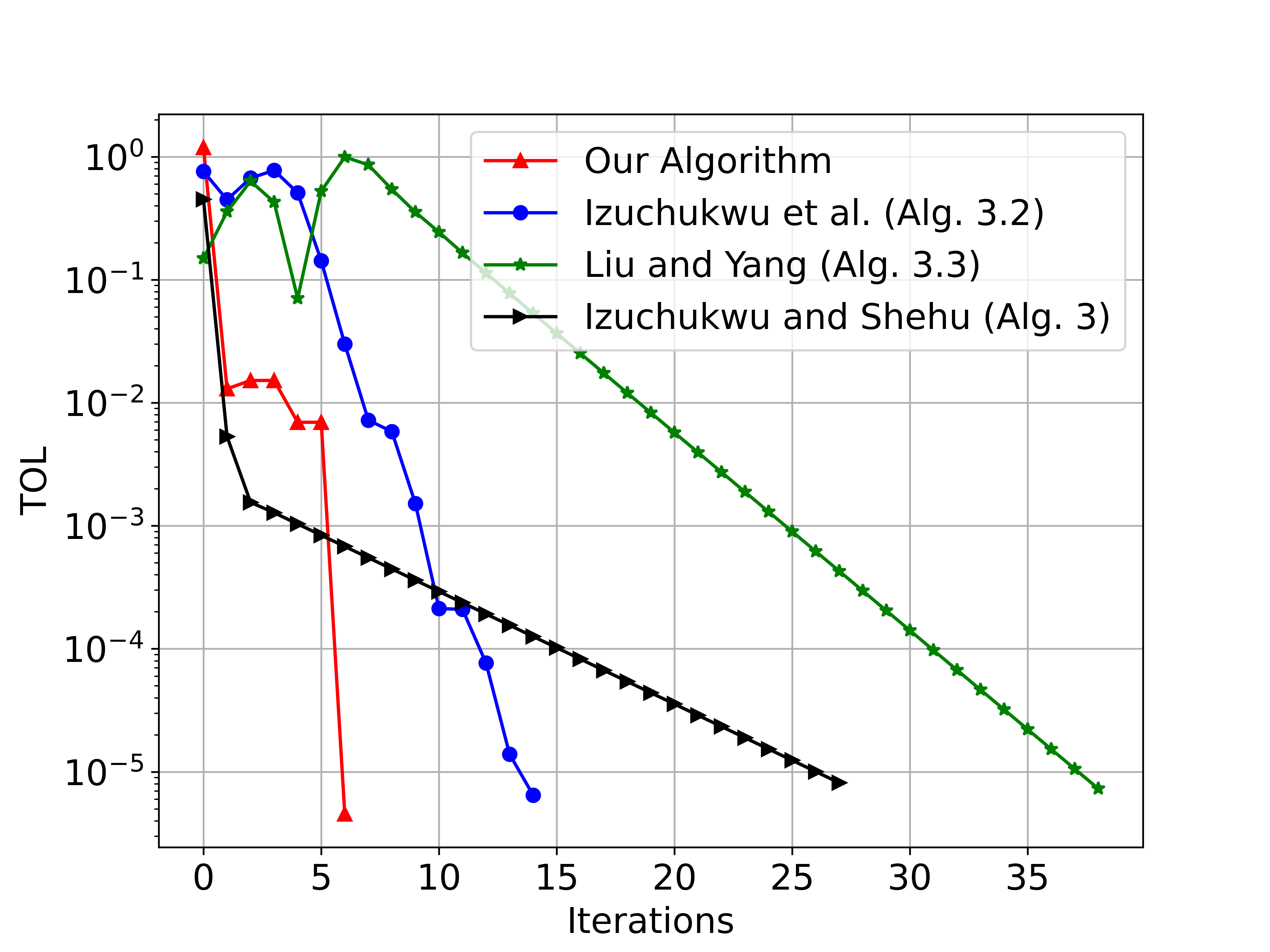}}}\quad
\subfigure{{\includegraphics[scale=0.35]{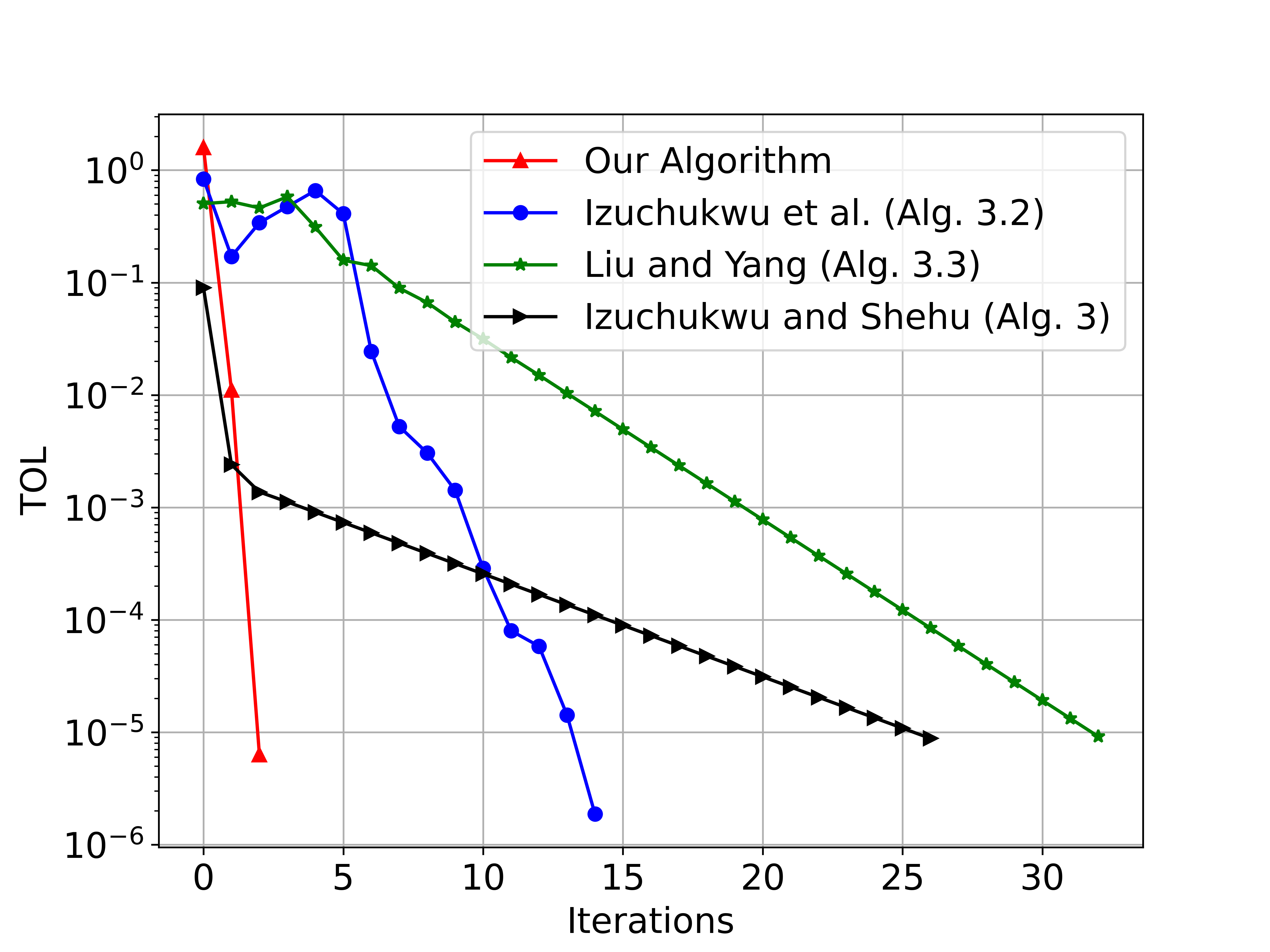}}}\quad
\subfigure{{\includegraphics[scale=0.35]{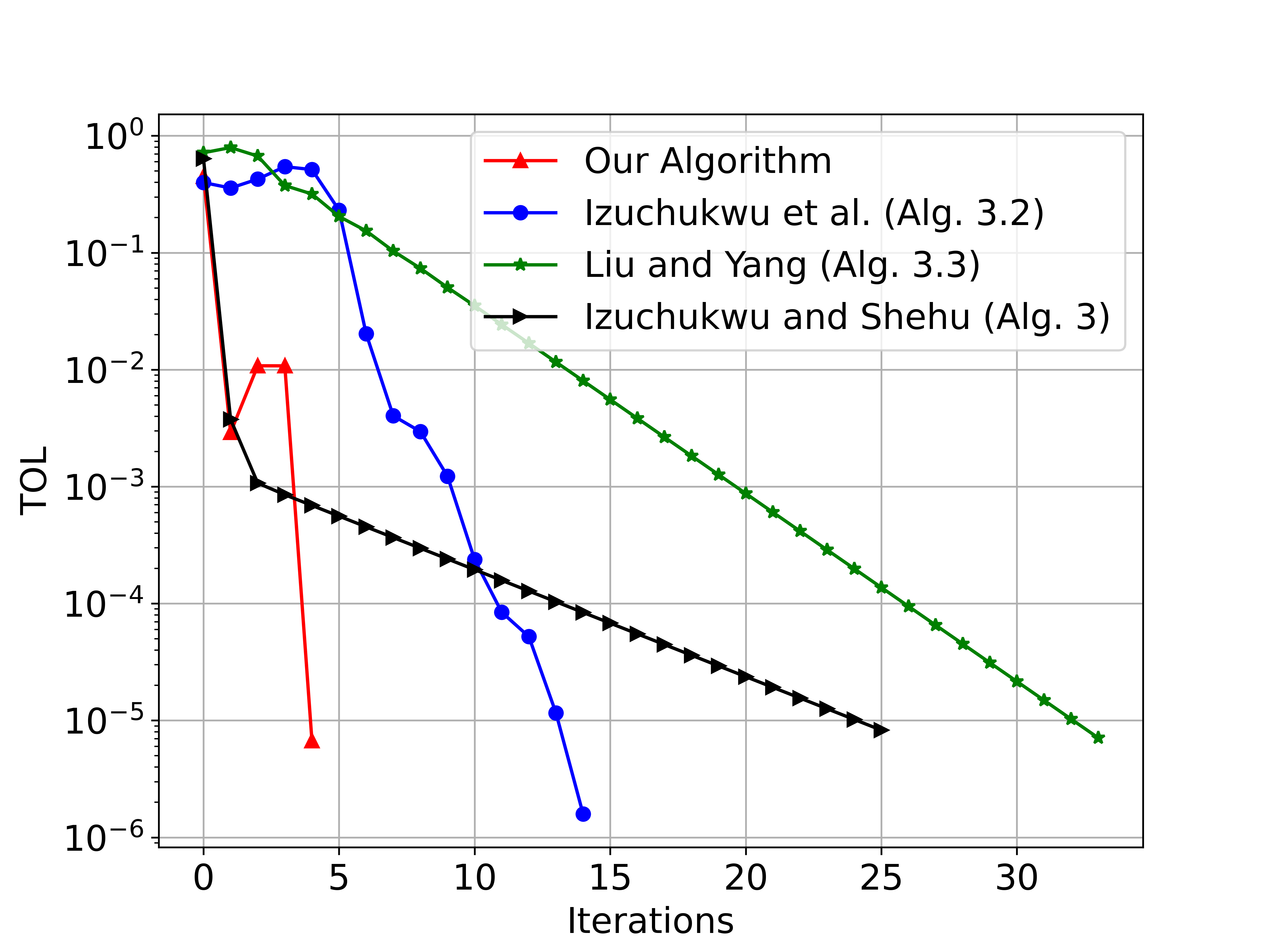}}}
\caption{Comparison of $\text{TOL}$ for Example \ref{exm_2} with $\epsilon=10^{-5}$: Top Left: Case 1; Top Right: Case 2; Bottom left: Case 3; Bottom right: Case 4}
 \label{fig:n=slr_ijcnn1}
\end{figure}

\begin{figure}[h]
\centering
\mbox{\subfigure{\includegraphics[scale=0.3]{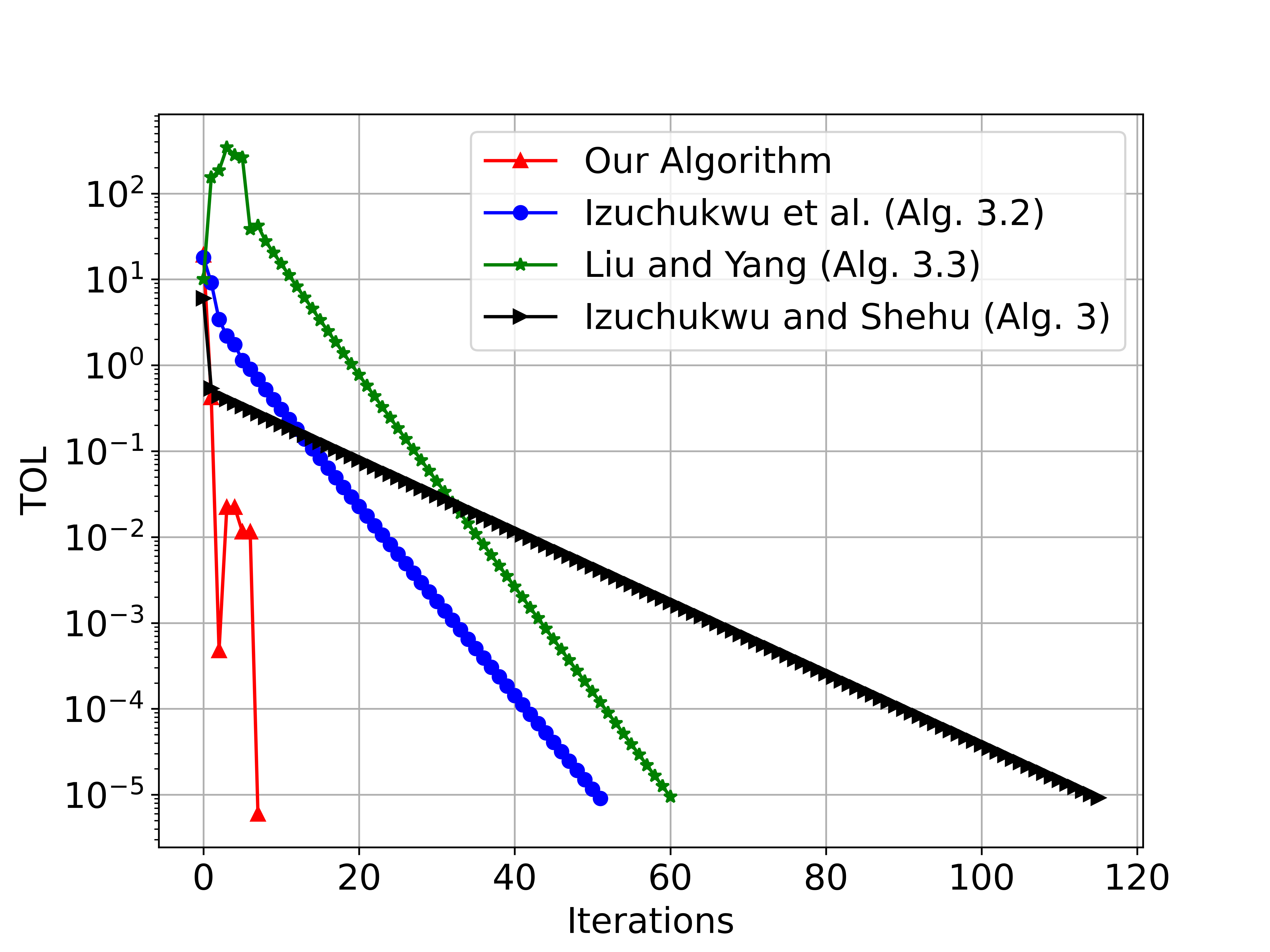} }\quad
\subfigure{{\includegraphics[scale=0.3]{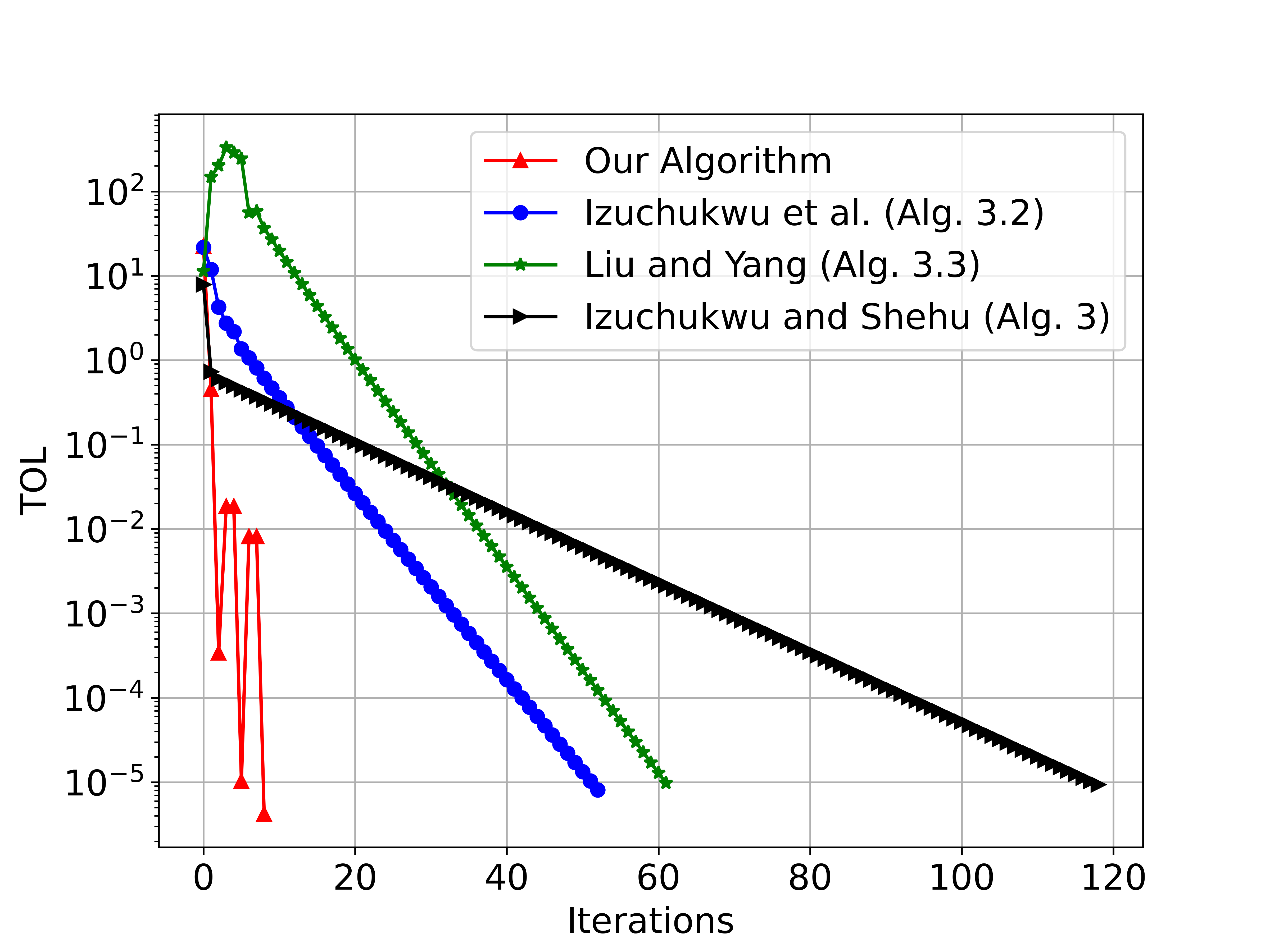} } }}\quad
{\subfigure{\includegraphics[scale=0.3]{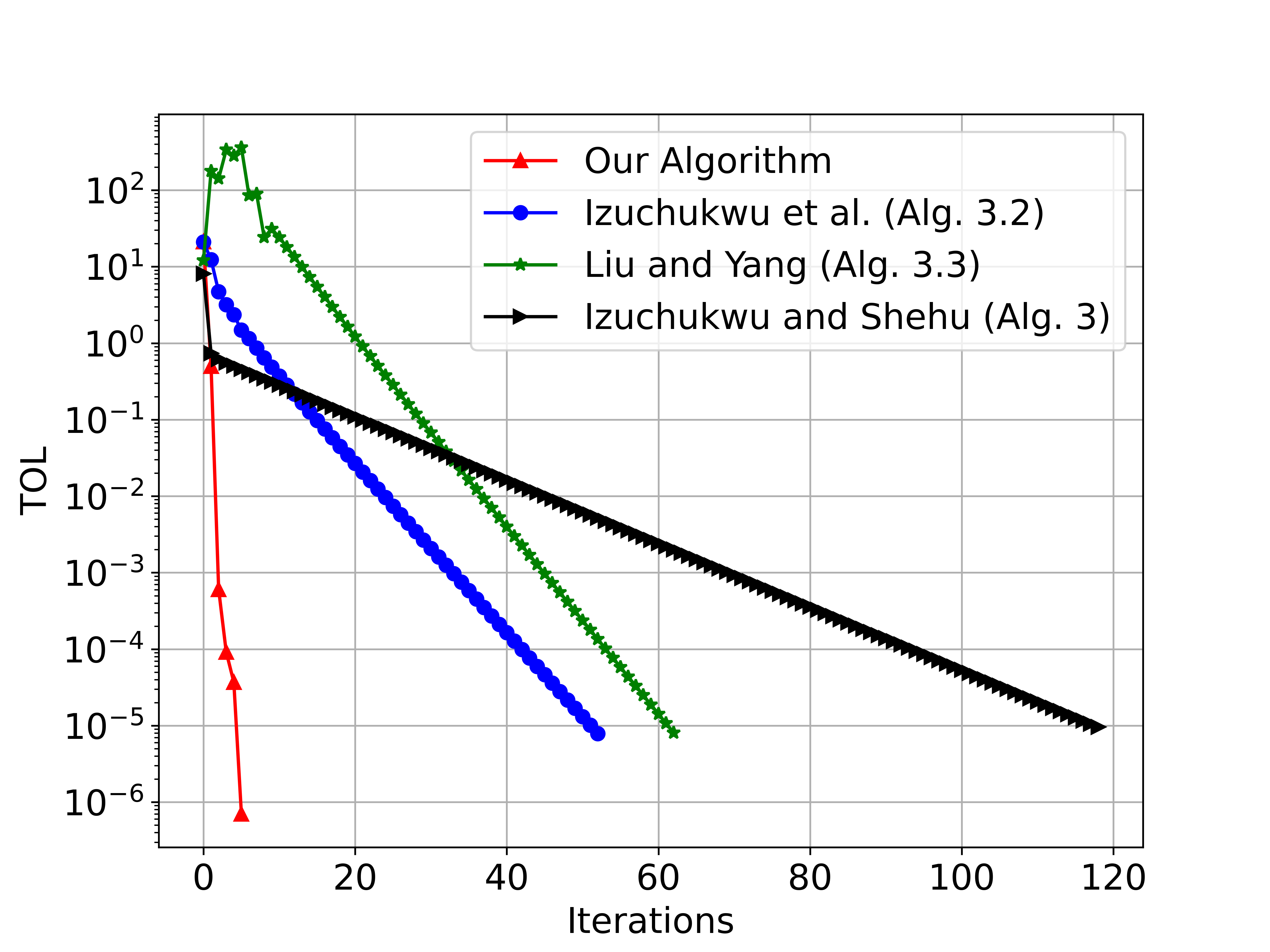} }\quad
\subfigure{{\includegraphics[scale=0.3]{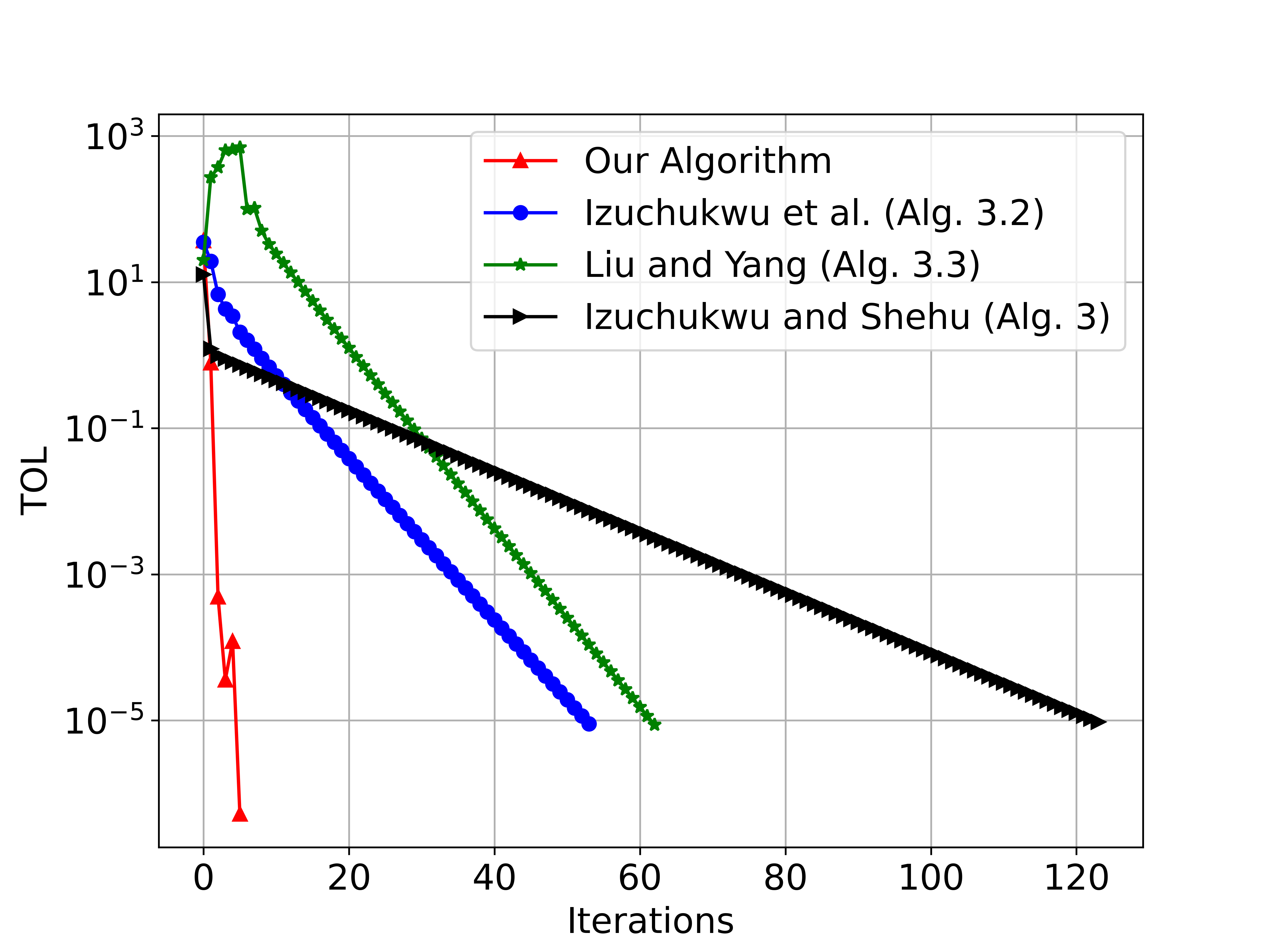} } }}
\caption{Comparison of $\text{TOL}$ for Example \ref{exm_3} with $\epsilon=10^{-5}$: Top Left: $m=50$; Top Right: $m=80$; Bottom left: $m=100$; Bottom right: $m=200$}
 \label{fig:n=slr_ijcnn1}
\end{figure}
\begin{figure}[h]
\centering
\mbox{\subfigure{\includegraphics[scale=0.3]{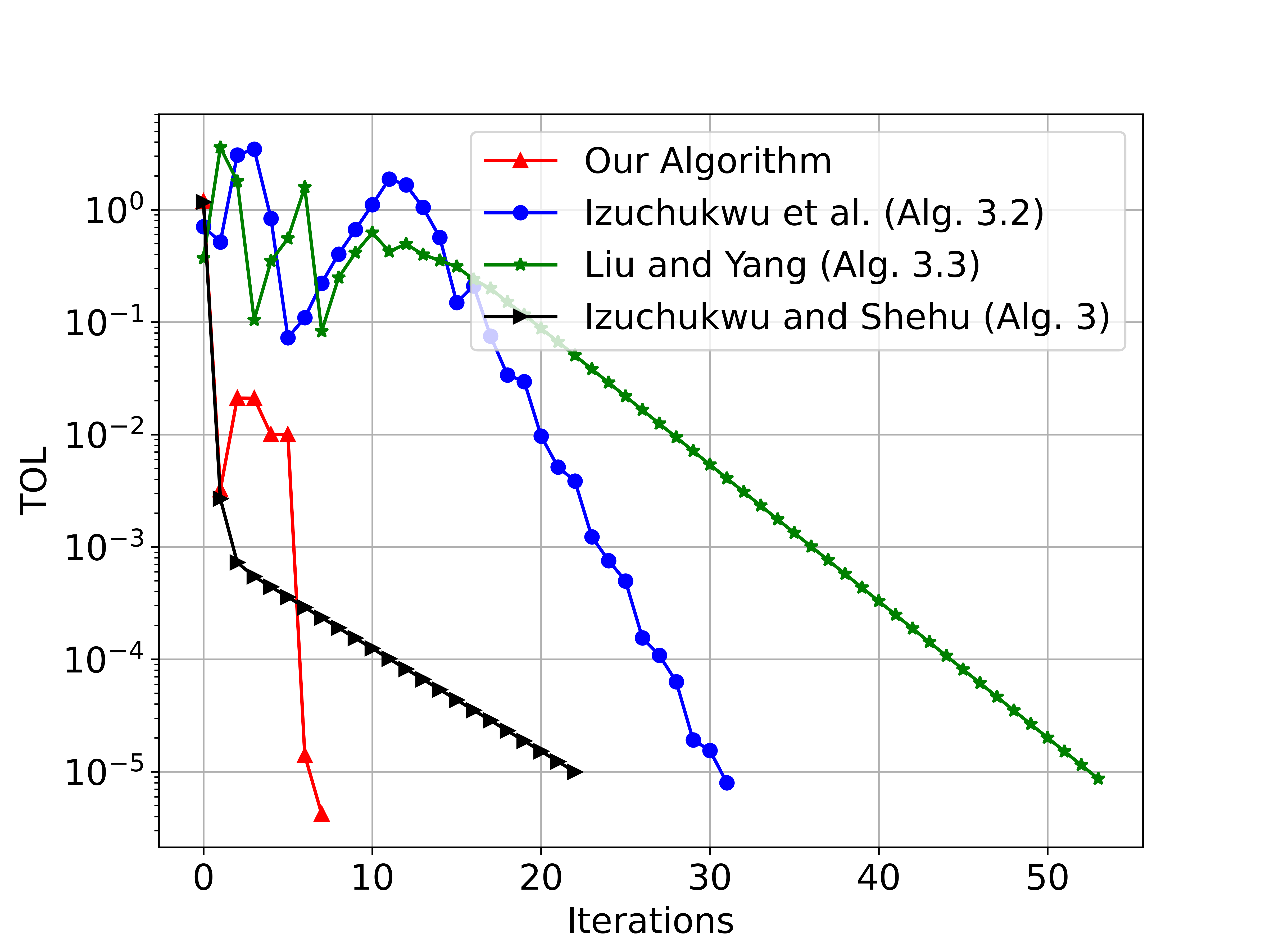} }\quad
\subfigure{{\includegraphics[scale=0.3]{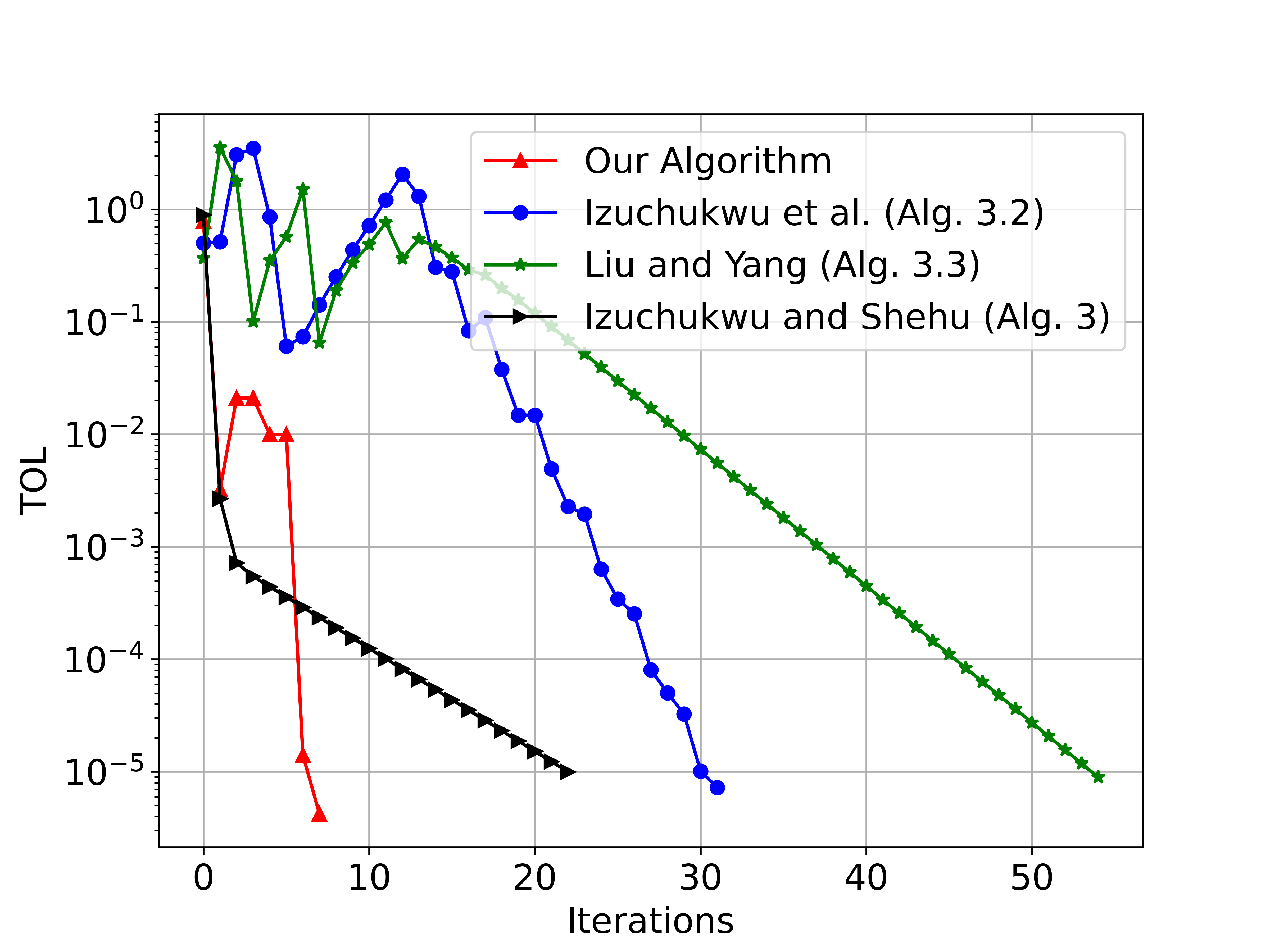} } }}\quad
{\subfigure{\includegraphics[scale=0.3]{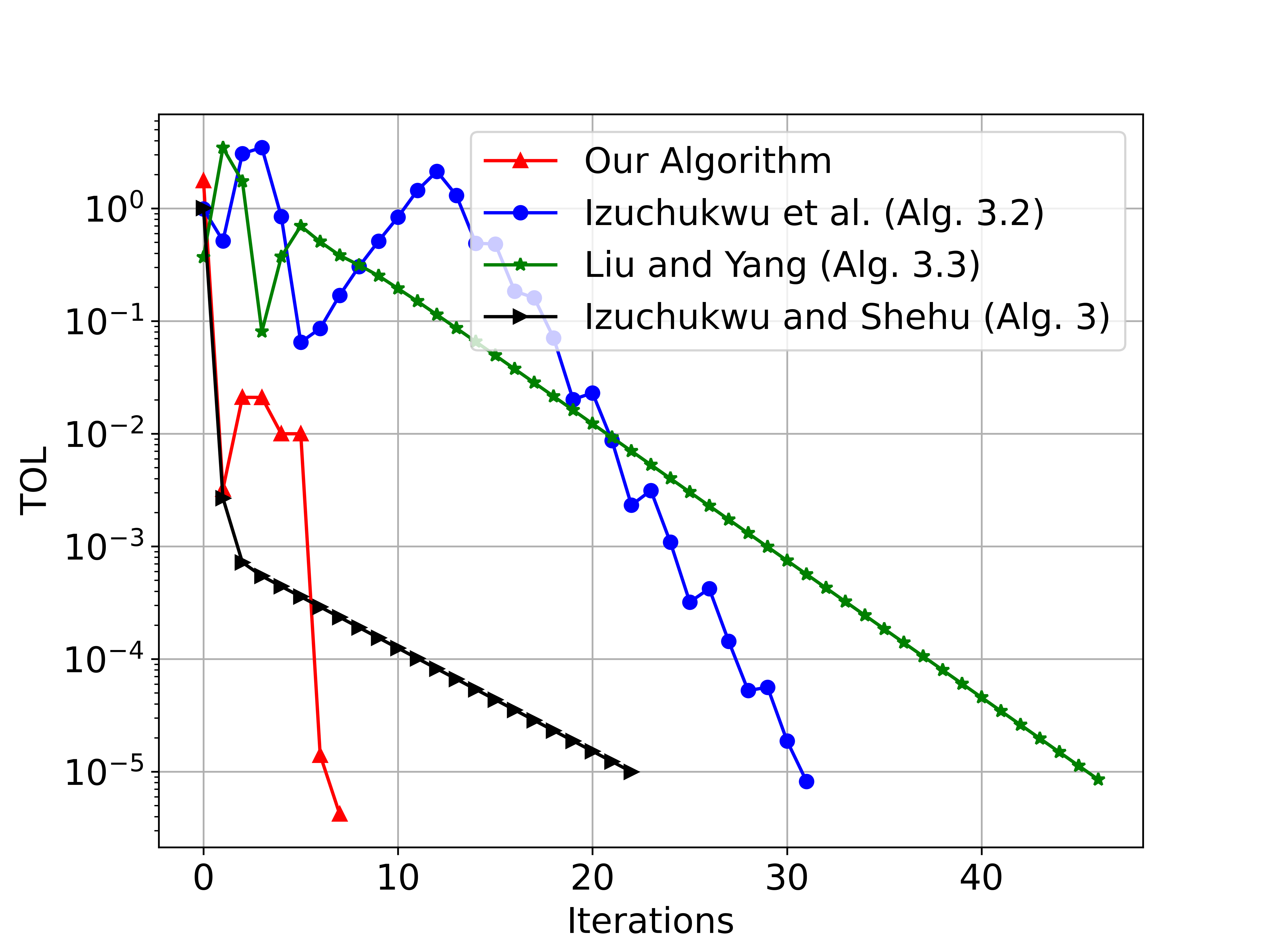} }\quad
\subfigure{{\includegraphics[scale=0.3]{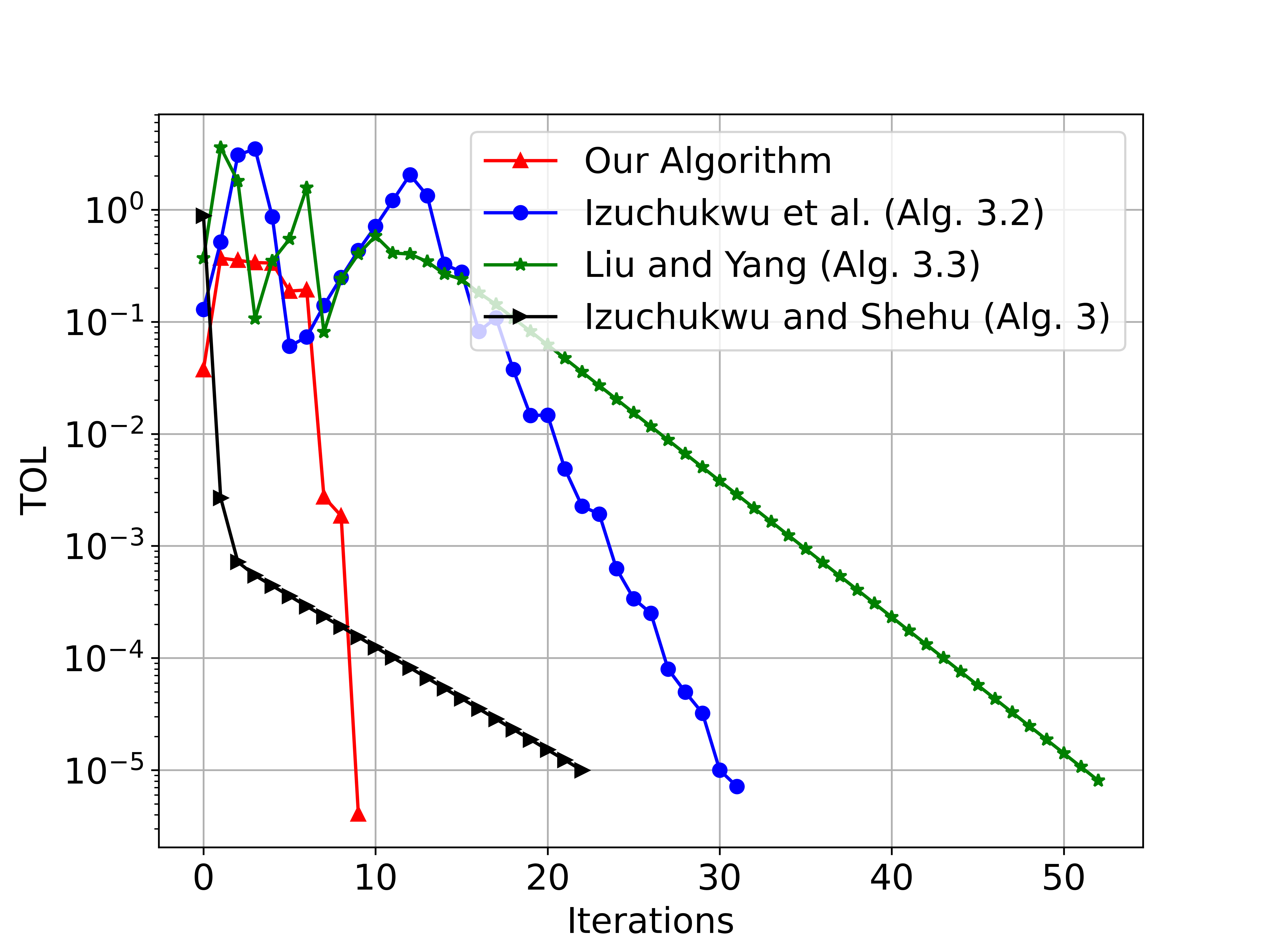} } }}
\caption{Comparison of $\text{TOL}$ for Example \ref{exm_4} with $\epsilon=10^{-5
}$: Top Left: Case 1; Top Right: Case 2; Bottom left: Case 3; Bottom right:Case 4}
 \label{fig:n=slr_ijcnn1}
\end{figure}
\begin{table}[htbp!]
\centering
\caption{Example \ref{exm_2}: Comparison of algorithms with $\epsilon = 10^{-5}$}
\begin{tabular}{lcccccccc}
\toprule
Algorithms & \multicolumn{2}{c}{Case 1} & \multicolumn{2}{c}{Case 2} & \multicolumn{2}{c}{Case 3} & \multicolumn{2}{c}{Case 4} \\
\cmidrule(r){2-3} \cmidrule(r){4-5} \cmidrule(r){6-7} \cmidrule(r){8-9}
 & CPU(s) & Iter & CPU(s) & Iter & CPU(s) & Iter & CPU(s) & Iter \\
\midrule
Our Algorithm & 0.0003 & 3 & 0.0010 & 7 & 0.0003 & 3 & 0.0006 & 5 \\
Algorithm 3 in \cite{izuchukwu2024golden}& 0.0012 & 27 & 0.0042 & 28 & 0.0013 & 28 & 0.0117 & 26 \\
Algorithm 3.2 in \cite{izuchukwu2023simple}& 0.0021 & 15 & 0.0016 & 35 & 0.0022 & 15 & 0.0143 & 15 \\
Algorithm 3.3 in \cite{liu2020weak}& 0.0040 & 29 & 0.0035 & 39 & 0.0033 & 33 & 0.0290 & 34\\
\bottomrule
\end{tabular}
\end{table}
\begin{table}[htbp!]
\centering
\caption{Example \ref{exm_4}: Comparison of algorithms with $\epsilon = 10^{-5}$}
\begin{tabular}{lcccccccc}
\toprule
Algorithms & \multicolumn{2}{c}{Case 1} & \multicolumn{2}{c}{Case 2} & \multicolumn{2}{c}{Case 3} & \multicolumn{2}{c}{Case 4} \\
\cmidrule(r){2-3} \cmidrule(r){4-5} \cmidrule(r){6-7} \cmidrule(r){8-9}
 & CPU(s) & Iter & CPU(s) & Iter & CPU(s) & Iter & CPU(s) & Iter \\
\midrule
Our Algorithm & 0.0102 & 8 & 0.0104 & 8 & 0.0025 & 8 & 0.0045 & 10 \\
Algorithm 3 in \cite{izuchukwu2024golden}& 0.0234 & 23 & 0.223 & 23 & 0.0052 & 23 & 0.0081 & 23 \\
Algorithm 3.2 in \cite{izuchukwu2023simple}& 0.0354 & 32 & 0.0291 & 32 & 0.0052 & 32 & 0.0108 & 32 \\
Algorithm 3.3 in \cite{liu2020weak}& 0.0691 & 54 & 0.0717 & 55 & 0.0176 & 47 & 0.0309 & 53\\
\bottomrule
\end{tabular}
\end{table}\clearpage
\subsection{Application to signal recovery}\noindent
Next, we perform the following signal processing application where our aim is to recover the true signal from a noisy observation. Consider the noisy signal $y$ given by
\begin{equation*}
    y = Bx + \mu,
\end{equation*}
where $B \in \mathbb{R}^{m \times n}$ is the observation matrix, $x \in \mathbb{R}^n $ represents the true sparse signal, and $\mu \in \mathbb{R}^m$ accounts for additive noise.
Nonetheless, to recover the true signal, we consider solving the following constrained Lasso optimization problem.
\begin{equation}\label{lasso}
    \min_{x\in{\mathbb{R}}^n}\frac{1}{2}\lVert Bx-y\rVert^2~\text{subject to}~\lVert x\rVert_1\leq l.
\end{equation}
Further, by using the first order optimality condition, (\ref{lasso}) can be written as the following variational inequality problem.
\[\text{Find}~p^*\in\mathscr{C}~\text{such that}~\langle \mathscr{A}p^*, q-p^*\rangle\geq 0~\forall q\in\mathscr{C},\]
where $\mathscr{A}p^*=B^\top(Bp^*-y)$ and the set $\mathscr{C}=\{x\in{\mathbb{R}}^n:\lVert x\rVert_1\leq l\}$.\par
\begin{figure}[htbp]
    \centering

    \subfigure[True signal]{
        \includegraphics[width=0.46\textwidth]{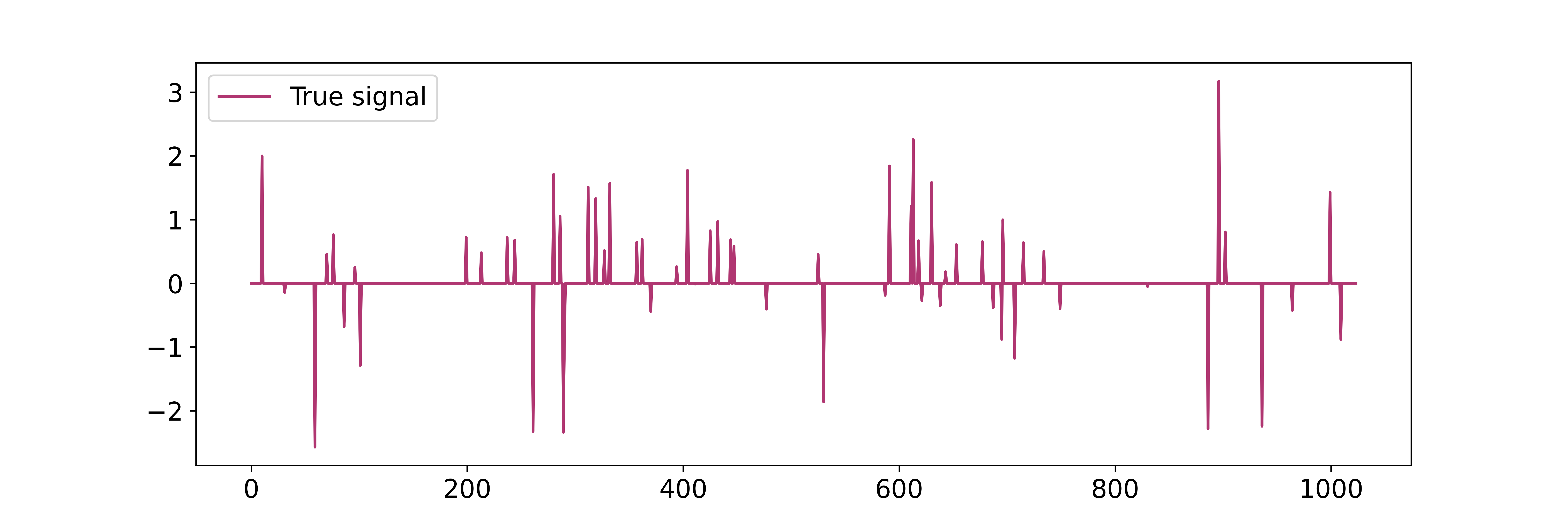}
        \label{fig:true}
    }
    \hfill
    \subfigure[Recovered by our method]{
        \includegraphics[width=0.46\textwidth]{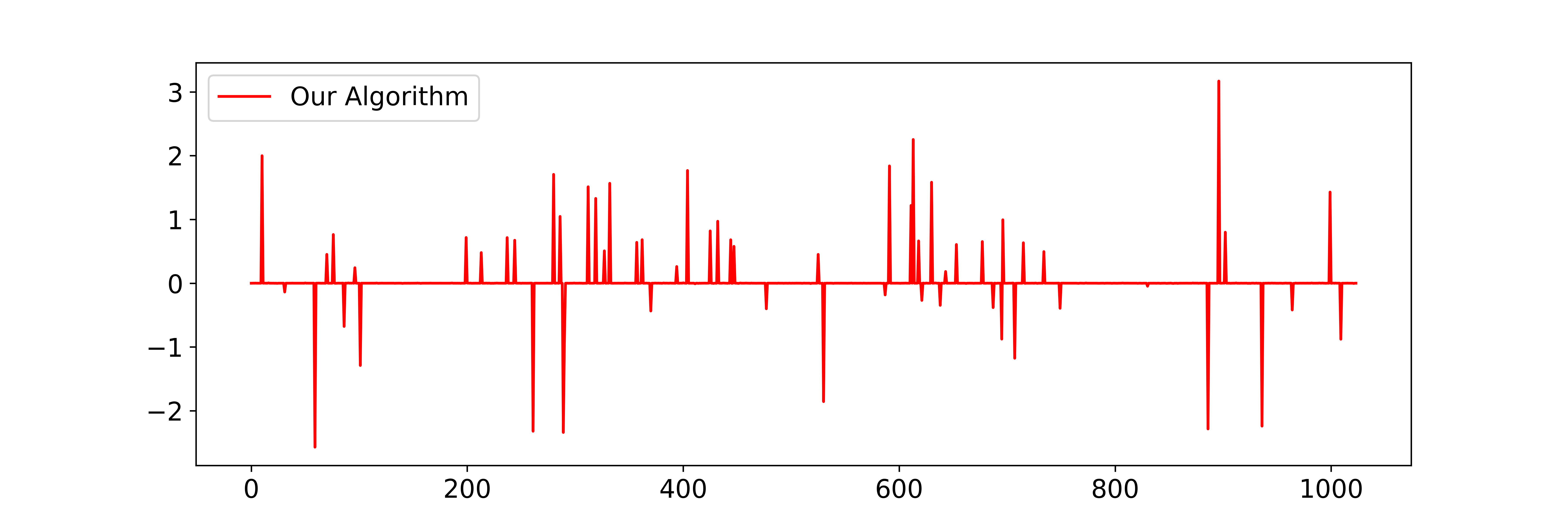}
        \label{fig:noisy}
    }
    \hfill
    \subfigure[Noisy signal]{
        \includegraphics[width=0.46\textwidth]{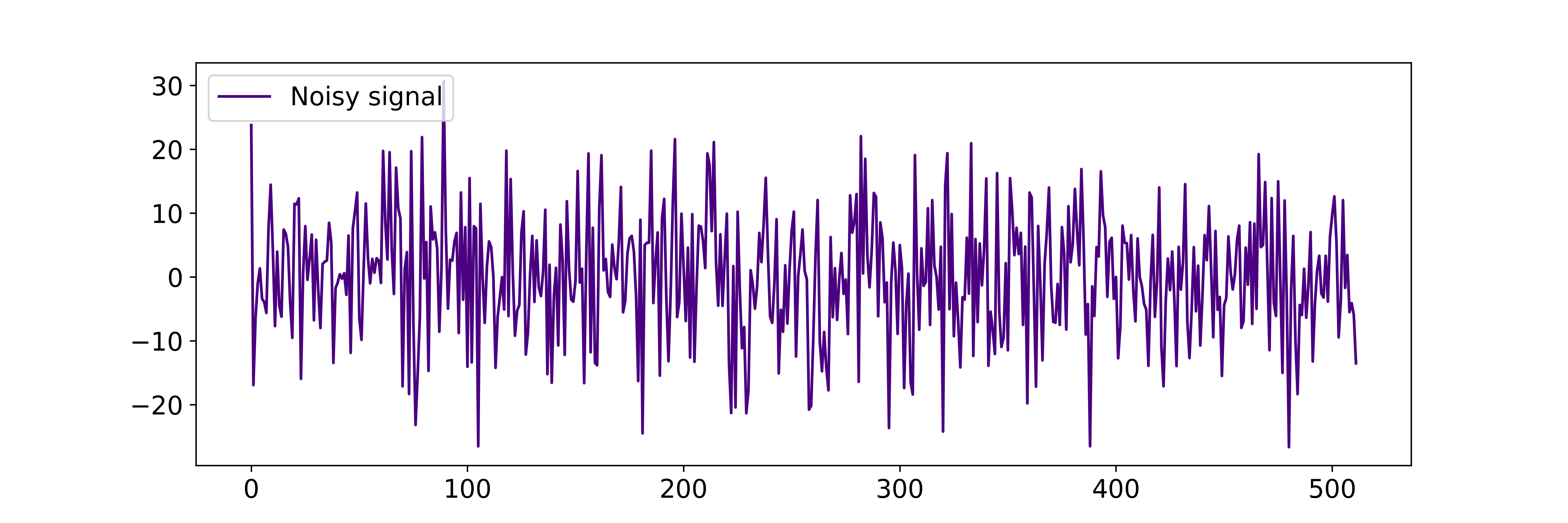}
        \label{fig:our}
    }

    \vspace{0.3cm}

    \subfigure[Recovered by Izuchukwu and Shehu]{
        \includegraphics[width=0.46\textwidth]{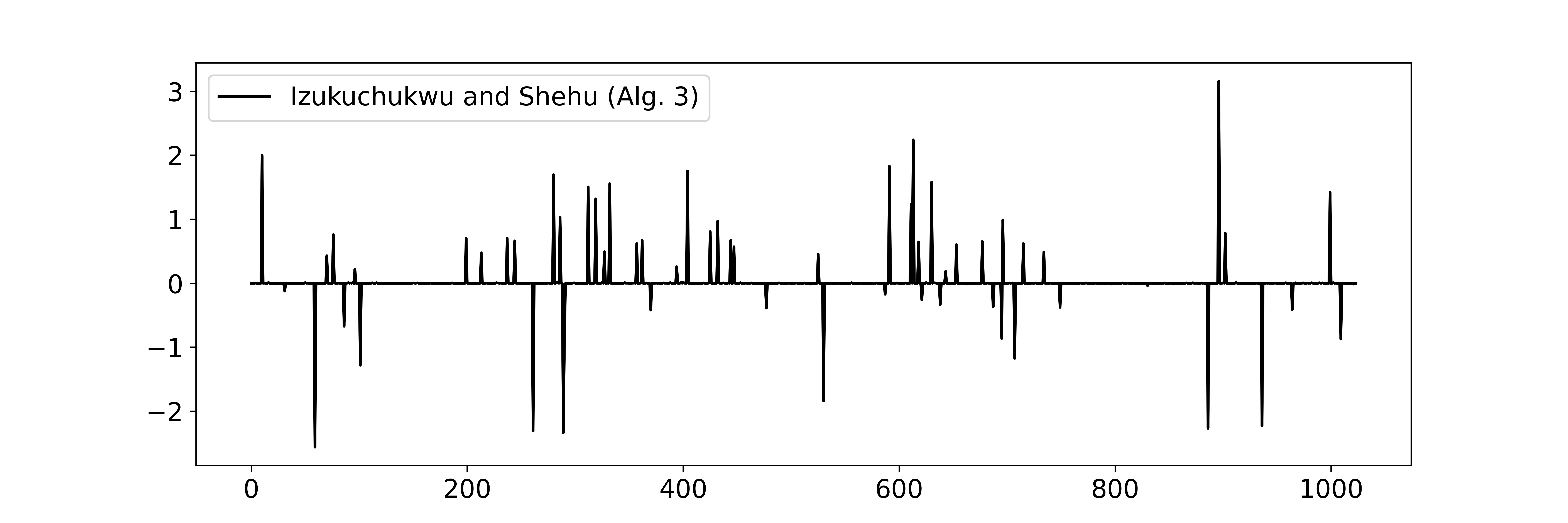}
        \label{fig:izu_shehu}
    }
    \hfill
    \subfigure[Recovered by Izuchukwu et al.]{
        \includegraphics[width=0.46\textwidth]{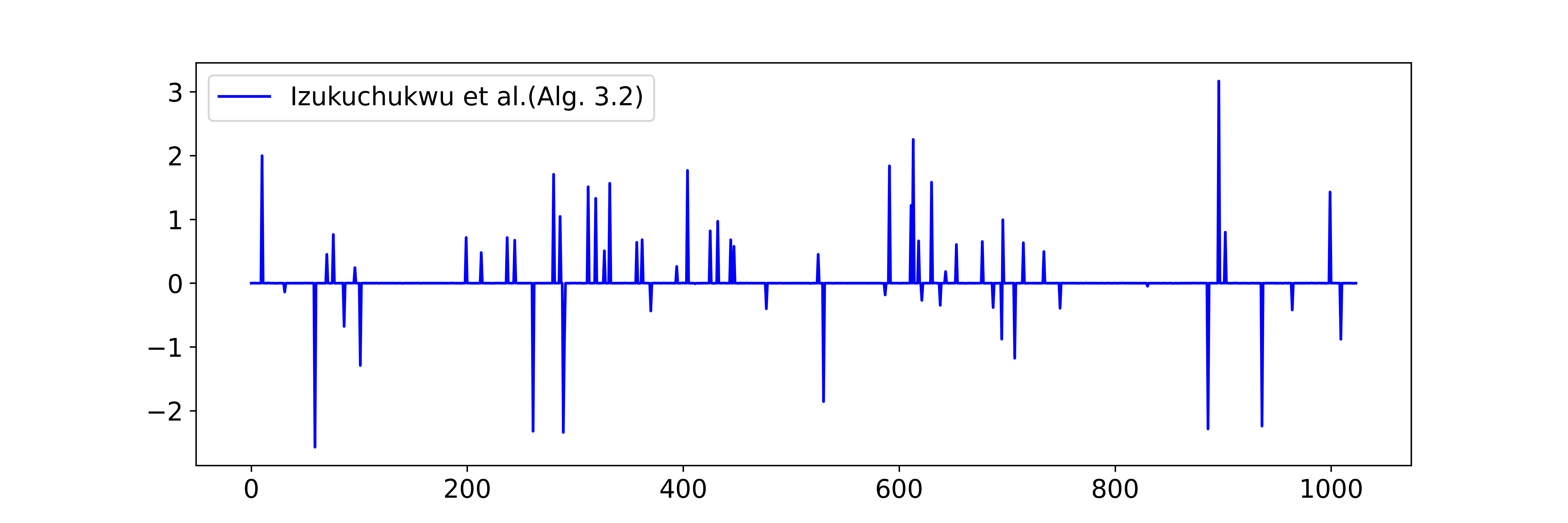}
        \label{fig:izu_etal}
    }

    \caption{Results of signal recovery by different algorithms}
    \label{fig:signal_comp}
\end{figure}
\begin{figure}[h!]
		\centering
		\includegraphics[scale=0.4]{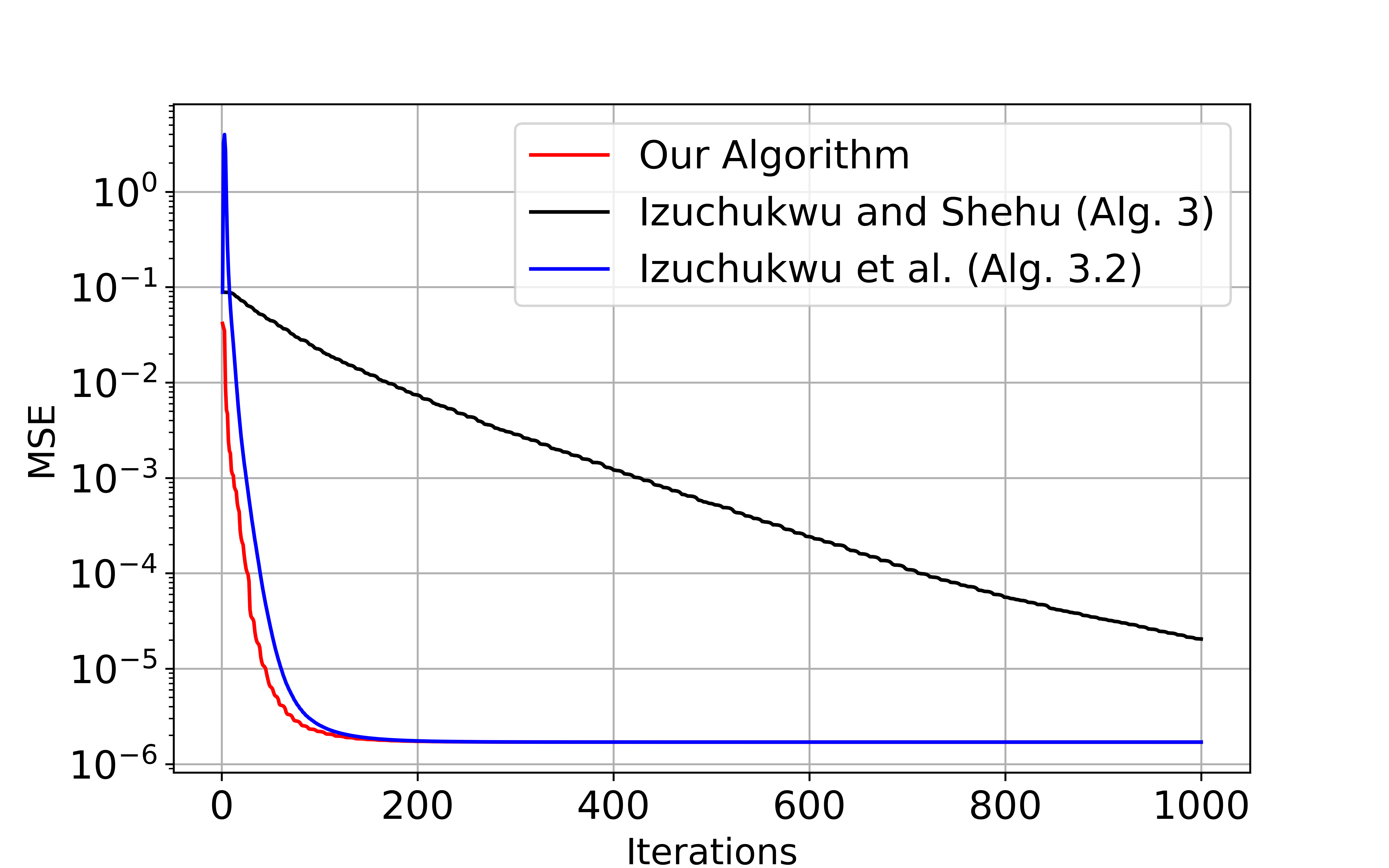}
		\caption{MSE comparison}\label{mse_comparison}
	\end{figure}

In our experiment, the matrix $B$ and the noise vector $\mu$ are generated randomly using \texttt{numpy library} in Python as: $B = \texttt{np.random.randn}(m, n)$ and $\mu = 10^{-3}~\texttt{np.random.randn}(m, 1)$. The true signal $x$ is constructed with $s$ non zero entries located at random positions.
To track the recovery accuracy, we use the mean square error (MSE) given by $\frac{1}{n} \|p^* - x\|^2$,
where \( p^* \) is the signal estimated by the algorithm.
All algorithms are initialized with $v_0 = \mathbf{0} \in \mathbb{R}^n$, and $v_1$ is initialized randomly. The iteration stops when the condition \(\text{MSE} < 10^{-6} \) is met or a maximum of $1000$ iterations is reached. In our simulations, we set $n = 1024, m = 512$, and take $l=s=60$. Figures \ref{fig:signal_comp} and \ref{mse_comparison} collect all the graphical results for this experiment.
\begin{table}[htbp!]
\centering
\caption{Example \ref{exm_3}: Comparison of algorithms with $\epsilon = 10^{-5}$}
\begin{tabular}{lcccccccc}
\toprule
Algorithms & \multicolumn{2}{c}{Case 1} & \multicolumn{2}{c}{Case 2} & \multicolumn{2}{c}{Case 3} & \multicolumn{2}{c}{Case 4} \\
\cmidrule(r){2-3} \cmidrule(r){4-5} \cmidrule(r){6-7} \cmidrule(r){8-9}
 & CPU(s) & Iter & CPU(s) & Iter & CPU(s) & Iter & CPU(s) & Iter \\
\midrule
Our Algorithm & 0.0049 & 8 & 0.0074 & 9 & 0.0058 & 6 & 0.0127 & 6 \\
Algorithm 3 in \cite{izuchukwu2024golden}& 0.0307 & 116 & 0.0524 & 119 & 0.0371 & 120 & 0.1011 & 124 \\
Algorithm 3.2 in \cite{izuchukwu2023simple}& 0.0244 & 52 & 0.0373 & 53 & 0.0482 & 53 & 0.0869 & 54 \\
Algorithm 3.3 in \cite{liu2020weak}& 0.0231 & 61 & 0.0288 & 62 & 0.0371 & 63 & 0.0563 & 63\\
\bottomrule
\end{tabular}
\end{table}
\section{Conclusion}\label{sec:5}\noindent
In this study, we proposed a modified projection algorithm incorporating two momentum terms and a nondecreasing step size strategy to solve variational inequality problems involving quasimonotone and Lipschitz continuous operators. We established weak convergence of the method under standard assumptions. Through a series of numerical experiments on four benchmark problems, we demonstrated that the proposed algorithm outperforms several existing methods in terms of efficiency and accuracy. Furthermore, we demonstrated its practical applicability by using it to solve a signal recovery problem. Future research direction includes on extending the applicability of the method to broader problem classes, such as mixed variational inequality problems and equilibrium problems.
\section*{Acknowledgements}\noindent
Gourav Kumar acknowledges the institute postdoctoral fellowship from IIT Madras (Order No. Acad/IPDF/R12/2025). Santanu Soe was supported by the Prime Minister’s Research Fellowship program, Ministry of Education, Goverment of India.

\section*{Code availability} The Python codes used for the numerical experiments are available upon request.

\end{document}